\documentstyle[amscd,amssymb,verbatim]{amsart}
\newcommand{\lrar}[1]{\begin{picture}(50,10)(-25,-5)                          
\put(-25,0){\vector(1,0){50}}
\put(0,5){\makebox(0,0)[b]{\mbox{$#1$}}}
\end{picture}}

\newcommand{\ldar}[1]{\begin{picture}(10,50)(-5,-25)
\put(0,25){\vector(0,-1){50}}
\put(5,0){\mbox{$#1$}}
\end{picture}}

\newcommand{\ldrar}[1]{\begin{picture}(50,50)(-25,-25)
\put(-25,25){\vector(1,-1){50}}
\put(5,0){\mbox{$#1$}}
\end{picture}}

\newcommand{\Hyp}{\operatorname{Hyp}}
\newcommand{\inv}{\operatorname{inv}}

\newcommand{\Tr}{\operatorname{Tr}}
\newcommand{\Nm}{\operatorname{Nm}}

\newcommand{\Gal}{\operatorname{Gal}}

\newcommand{\F}{{\Bbb F}}
\newcommand{\G}{{\Bbb G}}
\newcommand{\A}{{\Bbb A}}

\newcommand{\mmu}{{\Bbb \mu}}

\newcommand{\hra}{\hookrightarrow}
\newcommand{\lan}{\langle}
\newcommand{\ran}{\rangle}

\newcommand{\GG}{{\cal G}}
\newcommand{\CC}{{\cal C}}

\newcommand{\Frob}{\operatorname{Frob}}
\newcommand{\Div}{\operatorname{Div}}
\newcommand{\supp}{\operatorname{supp}}

\newcommand{\Proj}{\operatorname{Proj}}
\renewcommand{\P}{{\Bbb P}}

\newcommand{\si}{\sigma}

\newcommand{\de}{\delta}
\newcommand{\eps}{\epsilon}

\renewcommand{\ker}{\operatorname{ker}}
\newcommand{\im}{\operatorname{im}}
\newcommand{\D}{{\cal D}}
\numberwithin{equation}{section}

\newtheorem{thm}{Theorem}[subsection]
\newtheorem{prop}[thm]{Proposition}
\newtheorem{lem}[thm]{Lemma}
\newtheorem{cor}[thm]{Corollary}
\newenvironment{rem}{\vspace{3mm}\noindent
{\bf Remark.}}{\vspace{3mm}}

\newcommand{\Pf}{\noindent {\it Proof}}

\newcommand{\ov}{\overline}

\newcommand{\rk}{\operatorname{rk}}
\newcommand{\ra}{\rightarrow}

\newcommand{\FF}{{\cal F}}

\renewcommand{\a}{\alpha}
\renewcommand{\b}{\beta}

\newcommand{\la}{\lambda}

\newcommand{\C}{{\Bbb C}}

\newcommand{\Z}{{\Bbb Z}}
\newcommand{\Q}{{\Bbb Q}}

\newcommand{\Ga}{\Gamma}

\newcommand{\wt}{\widetilde}

\newcommand{\sub}{\subset}
\newcommand{\ed}{\qed\vspace{3mm}}

\title{Fourier transform over finite field and identities between Gauss sums}
\author{D. Kazhdan and A. Polishchuk}

\begin{document}
\maketitle

\bigskip

\centerline{\sc Introduction}

\bigskip

This paper is a continuation of \cite{EKP}. In \cite{EKP} 
we considered distributions of the form 
$$C\psi(Q(x))\prod_{j=1}^k\chi_j(P_j(x))$$
over local fields,
where $Q(x)$ is a rational function of several variables, 
$P_j$ are polynomials, $C\in\C^*$,
$\psi$ is a non-trivial additive character, $\chi_j$ are multiplicative
characters. The main problem we posed in \cite{EKP} is to determine
when the Fourier transform of such a distribution
is again a distribution of the same type. 
It turned out that there is a necessary condition: the map 
$x\mapsto dQ(x)$ (which we assume to be dominant) should be birational.
The simplest example when this condition is satisfied is the monomial case:
$Q(x_1,\ldots,x_k)=a\prod_{i=1}^k x_i^{n_i}$,
where all $n_i$ are non-zero integers, the degree of $Q$ is either $0$ or $2$.
In this case it is natural to consider $P_j=x_j$, so we are looking
at the distributions of the form
\begin{equation}\label{form}
C\psi(a\prod_{i=1}^kx_i^{n_i})\prod_{i=1}^k\chi_j(x_i).
\end{equation}
We have shown in \cite{EKP} that the Fourier transform
of this distribution is again of the same form if and only if
certain monomial identity between the values of gamma-functions is satisfied.
The structure of these identities is well-known. In the archimedian case
all such identities follow from the multiplication law and the functional
equation. In the non-archimedian case to get interesting identities 
we have to consider characters of various extensions of a given local field.
The monomial identites between gamma-functions of a collection characters are 
governed by the linear relations between the corresponding induced
representations of the Galois group (see \cite{De2}).

In this work we will consider an analogue of this picture over finite fields.
The intuition coming from representation theory tells us that in this case
instead of distributions we have to consider $l$-adique perverse sheaves.
Namely, the Goresky-MacPherson extension of perverse sheaves
allows us to associate irreducible perverse sheaves on $\A^n$ (equipped 
with an action of the Frobenius) to an expression
of the form (\ref{form}). We show that the analogue of the main result of 
\cite{EKP} holds in this context. For this we develop
the theory of monomial identities between gamma-functions of characters
over a finite field $\F_q$ (i.e. Gauss sums). More precisely,
we are interested in identities which hold universally over arbitrary
finite extension of our finite field $\F_q$.
In section \ref{gauss-sec} we show that such identities are governed by
linear relations between divisors (formal linear combinations of points) 
on the set of multiplicative characters of $\ov{\F_q}$ of finite order.
This reminds of the situation over $\C$ where gamma-function is meromorphic
and one has to look at its divisor of poles and zeroes. Over finite
field instead of poles and zeroes we have the jump in the absolute
value of gamma-function occuring at the trivial character. In order
to use this jump effectively we have to invoke the theory of
Kloosterman sheaves due to N. Katz (see \cite{K1}).

In section \ref{gormac-sec} we compute the traces and eigenvalues of
Frobenius acting on the stalks of perverse sheaves on $\A^n$ associated with
expressions of the form (\ref{form}). This allows us
to extract elementary identities for finite Fourier transform
from the identities involving geometric
Fourier trasforms of these perverse sheaves.
Here we observe another similarity with the situation over $\C$:
the conditions under which our perverse sheaves 
have zero stalks over the union of coordinate hyperplanes
are very similar to
the conditions garanteeing that the corresponding distributions over $\C$
have no poles (see proposition 4.8 of \cite{EKP}).
This phenomenon might suggest that there exists a generalized
Riemann-Hilbert correspondence which includes
sheaves and $D$-modules generated by exponents of rational functions.
One can consider as further evidence for existence of such a correspondence 
the work \cite{K2} where monodromy groups of ``exponential'' perverse sheaves
are compared with differential Galois groups, as well as the works
\cite{BE1} and \cite{BE2} where sheaf-theoretic methods of computing
determinants of cohomology are transferred into the realm of irregular 
connections.

Finally, in section \ref{norm-sec}
we study ``forms'' of the identities for the Fourier transform
of perverse sheaves considered above. More precisely, we can replace powers
of variables in the expression (\ref{form}) by norms in finite
extensions of our field. We show that all our results can be
generalized to this case.

\noindent
{\it Notation}. All our schemes are of finite type over a finite field
(or over its algebraic closure). By a sheaf on such a scheme $S$
we mean an object of the derived category of constructible
$l$-adique sheaves on $S$ defined in \cite{DeWeil}. If
$i:x\ra S$ is a closed point, $F$ is a sheaf on $S$ then
we denote the stalk $i^*F$ by $F_x$ or by $F|_x$.

\noindent
{\it Acknowledgment}. We would like to thank P. Etingof and
M. Kontsevich for helpful discussions. We are grateful to K. Vilonen
for providing us with a simplier proof of proposition \ref{pure}.
Our work was supported by NSF grants DMS-9700477 and DMS-9700458.

\section{Fourier transform}

\subsection{Definition}

Let $\F_q$ be the finite field with $q$ elements,
$\psi:\F_q\ra\C^*$ be a non-trivial additive character.
For every finite-dimensional vector space $V$ over $\F_q$
we define the Fourier transform of a function $f:V\ra\C$ by the
formula
$$\hat{f}(x^*)=\sum_{x\in V}f(x)\psi(\lan x^*,x\ran)$$
where $x^*\in V^*$, $\lan\cdot,\cdot\ran$ denotes the
natural pairing between $V^*$ and $V$.
We have the following formulas:
$$\hat{\hat{f}}(x)=q^nf(-x),$$
$$(\hat{f},\hat{g})=q^n(f,g)$$
where $n=\dim V$, the scalar product $(f,g)$ is defined by the formula
$$(f,g)=\sum_x f(x)\ov{g(x)}.$$

Let us choose a prime $l$ such that $(l,q)=1$ and
an identification $\ov{\Q}_l\simeq\C$. We
denote by $L_{\psi}$ the Artin-Schreier sheaf on $\A^1$ associated
with $\psi$. For every scheme $S$ and a morphism $f:S\ra\A^1$ we
denote $L_{\psi}(f)=f^*L_{\psi}$. Let 
$\D(S)$ denotes the derived category of constructible complexes 
of $l$-adic sheaves on $S$ (see \cite{DeWeil}).

We consider vector spaces $V$ and $V^*$ as schemes over $\F_q$. Then
the Fourier-Deligne transform $\FF=\FF_{\psi}:\D(V)\ra \D(V^*)$
is defined by the formula
$$\FF(A)=Rp_{V^*!}(p_V^*A\otimes L_{\psi}(\lan x^*,x\ran)[n])$$
where $p_V$ and $p_{V^*}$ are the projections of $V\times V^*$ onto
its factors.

One can associate with every
object $K\in\D(V)$ its trace function $t_K:V(\F_q)\ra\C$
by considering traces of the Frobenius acting on the fibers of $K$.
Then we have
$$t_{\FF(K)}=(-1)^n\hat{t}_K.$$

We refer to \cite{BBD} for the definition of perverse $l$-adique sheaves.
The important propery of the Fourier-Deligne transform is that it
sends perverse sheaves on $V$ to perverse sheaves on $V^*$ (see \cite{La2}).

\subsection{Gauss sums}

Let $\la:\F_q^*\ra\C^*$ be a non-trivial mulitplicative character
of $\F_q^*$. Then we can consider $\la$ as a function on $\F_q$
by setting $\la(0)=0$.
Then we have
$$\hat{\la}=g(\la)\la^{-1}$$
where 
$$g(\la)=g(\la,\psi)=\sum_{x\in\F_q^*}\la(x)\psi(x)$$
is the Gauss sum associated with $\la$.
Applying the Fourier transform twice we derive that
$$g(\la)g(\la^{-1})=\la(-1)q.$$
Also we have
$$\ov{g(\la)}=\la(-1)g(\la^{-1}).$$

It is convenient to extend this notation to the case of the trivial
character by setting
$$g(1)=\sum_{x\in\F_q^*}\psi(x)=-1.$$
Then for any multiplicative character $\la$ we have the following
formulas:
$$\sum_{x\in\F_q^*}\psi(xy)\la(x)=g(\la)\la^{-1}(y)$$
for any $y\in\F_q^*$,
$$g(\la^{-1})\equiv \la(-1)g(\la)^{-1}\mod q^{\Z}$$
where we denote by $q^{\Z}$ the subgroup in $\C^*$ consisting
integer powers of $q$.

Let $\la$ be a non-trivial character of $\F^*_q$.
There is a smooth sheaf $L_{\la}$ of rank 1 on $\G_m$ such that
$\la=t_{L_{\la}}$. If $j:\G_m\ra\A^1$ is the natural inclusion
then we have $Rj_*L_{\la}=j_!L_{\la}$. One has
$$\FF(j_!L_{\la}[1])\simeq G(\la)\otimes
j_!L_{\la^{-1}}[1]$$
where
\begin{equation}\label{Gla}
G(\la)=G(\la,\psi)=H^1_c(\G_m, L_{\la}\otimes L_{\psi})
\end{equation}
is a one-dimensional $\ov{\Q}_l$-space equipped with an
action of $\Gal(\ov{\F}_q/\F_q)$ such that $\Frob_q$ acts
as $-g(\la,\psi)$. It is convenient to extend the notation to
the case of the trivial character by setting
\begin{equation}\label{G1}
G(1)=G(1,\psi)=H^1_c(\G_m, L_{\psi})=\ov{\Q}_l
\end{equation}
(so the Galois action on this space is trivial).

Consider the extension of fields $\F_q\sub\F_{q^d}$ of degree $d$.
Then we have 
$$L_{\psi}\otimes_{\F_q}\F_{q^d}\simeq L_{\psi\circ \Tr_d}$$
$$L_{\la}\otimes_{\F_q}\F_{q^d}\simeq L_{\la\circ \Nm_d}$$
where $\Tr_d:\F_{q^d}\ra\F_q$ is the trace,
$\Nm_d:\F_{q^d}^*\ra \F_q^*$ is the norm homomorphism.
It follows that $G(\la\circ \Nm_d,\psi\circ\Tr_d)$ as 
$\Gal(\ov{\F}_q/\F_{q^d})$-representation
is obtained from $\Gal(\ov{\F}_q/\F_q)$-representation
$G(\la,\psi)$ by restriction. Thus,
we arrive to the identity
\begin{equation}\label{HasDavsim}
-g(\la\circ \Nm_d,\psi\circ\Tr_d)=(-g(\la,\psi))^d
\end{equation}
which is due to Hasse and Davenport (see \cite{DH}).

Another important identity for Gauss sums is the Hasse-Davenport product 
formula (see \cite{DH}):
\begin{equation}\label{HasDav}
\frac{g(\la^n)}{g(1)}=
\la(n^n)\cdot\prod_{i=0}^{n-1}\frac{g(\la\eps_n^i)}{g(\eps_n^i)}
\end{equation}
where $\eps_n$ is a character of order $n$ on $\F^*_q$
(so $n|(q-1)$, in particular, $n\neq 0$ in $\F_q$). 
The geometric interpretation of this identity is more involved
(see \cite{K1}, prop. 5.6.2) and it seems unlikely that
it admits an elegant geometric proof.  

\subsection{Main lemma}

Let us consider the standard embedding $j:\G_m^n\ra\A^n$. 
Assume that we have simple perverse sheaves $K$ and $K'$ on $\G_m^n$.
We are looking for a criterion checking that
$$\FF(j_{!*}K)\simeq j_{!*}K'.$$

For every $q_1=q^d$ let us denote by $f_{q_1}$ (resp. $f'_{q_1}$)
the trace function of $K\otimes \F_{q_1}$ (resp. $K'\otimes \F_{q_1}$).
These are $\C$-valued functions on $(\F_{q_1}^*)^n$. 
For every collection of characters $\la_1,\ldots,\la_n$
of $\F_{q_1}^*$ let us denote by $\la_1\otimes\ldots\otimes\la_n$ 
the corresponding
function $\la_1(x_1)\ldots\la_n(x_n)$ on $(\F_{q_1}^*)^n$.

\begin{lem}\label{mainlemma} 
Assume that for every $d>0$ and for every collection of
non-trivial characters $(\la_1,\ldots,\la_n)$ of
$\F_{q_1}^*$, where $q_1=q^d$ one has
$$(-q_1)^n(f_{q_1},\la_1\otimes\ldots\otimes\la_n)=
\prod_{i=1}^n \ov{g(\la_i,\psi\circ\Tr_d)} \cdot 
(f'_{q_1},\la_1^{-1}\otimes\ldots
\otimes\la_n^{-1})$$
and that this number is not zero for at least one collection
of non-trivial characters $(\la_1,\ldots,\la_n)$. Then
$$\FF(j_{!*}K)\simeq j_{!*}K'.$$
\end{lem}

\Pf . Since $K$ and $K'$ are irreducible it suffices to prove that
$j^*\FF(j_{!*}K)\simeq K'$.
For every scheme $S$ let us denote by $K_0(S)$ the Grothendieck group
of the category of perverse sheaves on $S$. The Fourier-Deligne transform
induces a homomorphism $\FF:K_0(\A^n)\ra K_0(\A^n)$. 
Let us denote by $\ov{K_0(\G_m^n)}$ the quotient of $K_0(\G_m^n)$
by the sum of the subgroups $p_i^*(K_0(\G_m^{n-1}))$ where 
$p_i:\G_m^n\ra\G_m^{n-1}$ is the projection omitting the $i$-th factor.
Then $\FF$ induces a well defined homomorphism
$$\ov{\FF}:\ov{K_0(\G_m^n)}\ra\ov{K_0(\G_m^n)}.$$
Indeed, for a sheaf $A$ on $\G_m^n$ we can choose a sheaf
$\wt{A}$ on $\A^n$ such that $A\simeq j^*\wt{A}$ and set
$\ov{\FF}([A])=[j^*\FF(\wt{A})]$. Since $K_0(\G_m^n)$ is the
quotient of $K_0(\A^n)$ by the subgroup generated by sheaves
supported on coordinate hyperplanes and the Fourier transform
interchanges such sheaves with sheaves of the form $p_i^*A$, the
map $\ov{\FF}$ is well-defined. The involutivity of the Fourier
transform implies that $\ov{\FF}$ is an isomorphism.
We have a natural basis in $\ov{K_0(\G_m^n)}$ corresponding to
simple perverse sheaves on $\G_m^n$ 
which are not constant on any factor $\G_m$
(i.e. do not belong to $p_i^*\D(\G_m^{n-1})$ for any $i$), and
the map $\ov{\FF}$ induces some permutation on this basis.

Our assumption that the scalar product of $f_{q_1}$ (resp. $f'_{q_1}$)
with some non-trivial multiplicative characters is non-zero implies
that $K$ (resp. $K'$) is not constant on any factor $\G_m$. Hence,
$\ov{\FF}([K])$ and $[K']$ are both elements of the basis in 
$\ov{K_0(\G_m^n)}$. 

Note that a function $f$ on $(\F_q^*)^n$ is completely 
determined by the scalar products 
with functions of the form $\la_1\otimes\ldots\otimes \la_n$,
where $\la_i$ are characters of $\F_q^*$.
Moreover, $f$ can be represented in the form 
$$f=\sum_{i=1}^n p_i^*f_i$$
if and only if 
$$(f,\la_1\otimes\ldots\otimes \la_n)=0$$
for all collections of non-trivial characters $\la_1,\ldots,\la_n$.
Now we claim that an element $x\in K_0(\G_m^n)$ 
lies in $\sum_i p_i^*K_0(\G_m^{n-1})$ if and only if
for all extensions $\F_q\subset\F_{q_1}$ the trace function
$t_x$ of $x$ over $\F_{q_1}$ has form $\sum_i p_i^* f_i$.
Indeed, a function $f$ can be represented in the form $\sum_i p_i^*f_i$
if and only if the following equality holds:
$$\sum_{I\sub [1,\ldots,n]}(-1)^{|I|} p_I^*\si_I^*f=0$$
where $p_I:\G_m^n\ra\G_m^I$ is the natural projection,
$\si_I:\G_m^I\ra\G_m^n$ is the embedding which sends $(a_i)_{i\in I}$ to
$(b_i)_{1\le i\le n}$ where $b_i=a_i$ for $i\in I$, $b_i=1$ otherwise.
Thus, our condition implies that all the trace functions of the element
$$\sum_{I\sub [1,\ldots,n]}(-1)^{|I|} p_I^*\si_I^*x$$
are zero. Hence, this element is zero by Theorem 1.1.2 of \cite{La2}.
 
Thus, an element of
$\ov{K_0(\G_m^n)}$ is completely determined by the scalar
products of its trace functions with $\la_1\otimes\ldots\otimes \la_n$,
where $\la_i$ are non-trivial characters. 

Now if $x\in \ov{K_0(\G_m^n)}$ and $(\la_1,\ldots,\la_n)$
is a collection of non-trivial characters of $\F_q^*$
then
$$(-q)^n(t_x,\la_1\otimes\ldots\otimes\la_n)=
\prod_i \ov{g(\la_i)}(t_{\ov{\FF}(x)},
\la_1^{-1}\otimes\ldots\otimes\la_n^{-1}).
$$
Thus, our assumptions imply that
$\ov{\FF}([K])=[K']$ in $\ov{K_0(\G_m^n)}$, hence
$j^*(\FF(j_{!*}K))\simeq K'$.
\ed

\section{Computation of Goresky-MacPherson extensions}
\label{gormac-sec}

\subsection{Vanishing of stalks}

Let $\psi$
be a non-trivial additive character of $\F_q$, $\chi_1,\ldots,\chi_k$
are characters of $\F_q^*$, $n_1,\ldots,n_k$ are
integers such that $(n_i,q)=1$. Then we define
$F^{n_1,\ldots,n_k}_{\chi_1,\ldots,\chi_k}$ as the Goresky-MacPherson
extension to $\A^k$ of the smooth perverse sheaf
$$L(\psi(\prod_i x_i^{n_i})\prod_i\chi_i(x_i))[k]:=
L_{\psi}(\prod_i x_i^{n_i})\otimes 
L_{\chi_1}\boxtimes\ldots\boxtimes L_{\chi_k}[k]$$
on $\G_m^k$. 
We are interested in computing the stalk 
$(F^{n_1,\ldots,n_k}_{\chi_1,\ldots,\chi_k})_0$ where
$0=(0,\ldots,0)\in\A^k$.
Of course, the interesting case is when some of $n_i$ are negative.

Let us also denote
$$\wt{F}^{n_1,\ldots,n_k}_{\chi_1,\ldots,\chi_k}
=Rj_*(L(\psi(\prod_i x_i^{n_i})\prod_i\chi_i(x_i)))$$
where $j$ is the embedding into $\A^k$ of the complement to the
subspace $x_i=0$ for all $i$ such that $n_i<0$ or $\chi_i\neq 1$,
the sheaf $L(\psi(\prod_i x_i^{n_i})\prod_i\chi_i(x_i)))$ is smooth
on this complement.

Consider the restriction of $\wt{F}^{n_1,\ldots,n_k}_{\chi_1,\ldots,\chi_k}$
to the open subset $x_1\neq 0$. Passing to the \'etale covering
$\ov{x_1}=x_1^{1/n_2}$ we can make a change of variables
$x_2'=x_2\ov{x_1}^{n_1}$, $x'_i=x_i$ for $i>2$, so that
$$\psi(\prod_i x_i^{n_i})\prod_i\chi_i(x_i)=
(\chi_1^{n_2}\chi_2^{-n_1})(\ov{x_1})\cdot
\psi(\prod_{i\ge 2} (x'_i)^{n_i})\prod_{i\ge 2}\chi_i(x'_i).$$
Thus, the pull-back of our sheaf under this \'etale covering
is the external tensor product of $L_{\chi_1^{n_2}\chi_2^{-n_1}}$
and $\wt{F}^{n_2,\ldots,n_k}_{\chi_2,\ldots,\chi_k}$.
Therefore, the calculation of all the stalks of
$\wt{F}^{n_1,\ldots,n_k}_{\chi_1,\ldots,\chi_k}$
reduces to the calculation of the stalks
$(\wt{F}^{n_{i_1},\ldots,n_{i_l}}_{\chi_{i_1},\ldots,\chi_{i_l}})_0$
for all subsets $\{i_1,\ldots,i_l\}\sub\{1,\ldots,k\}$.
Similar reasoning works for the sheaves
$F^{n_1,\ldots,n_k}_{\chi_1,\ldots,\chi_k}$.

In the case $k=1$ we have
$$(F^{-n}_{\chi})_0=(\wt{F}^{-n}_{\chi})_0=0$$
for any $n>0$ and any character $\chi$. This follows from the fact that
the Swan conductor of $L_{\psi}(1/x^n)$ at $x=0$ is equal to $n$.
Now let us consider the case $k=2$. Part of the following theorem
is contained in Theorem 3.1.1 of \cite{La1} (which is due to Deligne).
Our statement is slightly more general, however,
the proof is based on the same idea.

\begin{thm}\label{2dim}
Let $m,n>0$. Set $m=m'd$, $n=n'd$, where $d=gcd(m,n)$. Then

\noindent
(i) $(\wt{F}^{-m,-n}_{\chi_1,\chi_2})_0=
(F^{-m,-n}_{\chi_1,\chi_2})_0=0$ for any $\chi_1,\chi_2$,

\noindent
(ii) $(\wt{F}^{m,-n}_{\chi_1,\chi_2})_0=
(F^{m,-n}_{\chi_1,\chi_2})_0=0$ if $\chi_1^{n'}\chi_2^{m'}\neq 1$,

\noindent
(iii) $(F^{m,-n}_{\chi^{m'},\chi^{-n'}})_0=
H^1_c(\A^1,L_{\psi}(t^d)\otimes L_{\chi})[1]$,
where $L_{\chi}$ is extended to $\A^1$ by zero for $\chi\neq 1$
while $L_1=\ov{\Q}_{l,\A^1}$. If $\chi=1$ the dimension
of this cohomology space is $d-1$, while for $\chi\neq 1$ it is equal
to $d$.
\end{thm}

\Pf . Consider the Galois covering 
$$r:\A^2\ra\A^2:(\ov{x},\ov{y})\mapsto (x=\ov{x}^{n'},
y=\ov{y}^{m'})$$ 
with Galois group $G=\mmu_{n'}\times\mmu_{m'}=\mmu_{m'n'}$.
Then we have
$$(\wt{F}^{\pm m,-n}_{\chi_1,\chi_2})_0=
H^0(G,(\wt{F}^{\pm m'n'd,m'n'd}_{\chi_1^{n'},\chi_2^{m'}})_0).$$
This shows that for the proof of (i) and (ii) it suffices to
consider the case $m=n$.

\noindent
(i) Let $\pi:\wt{\A^2}\ra\A^2$ be the blow-up of $\A^2$ at $0$. 
Then there is a natural embedding
$\wt{j}:\A^1\times\G_m\ra\wt{\A^2}$ and we have
\begin{equation}\label{rga}
(\wt{F}^{0,-n}_{\frac{\chi_1}{\chi_2},\chi_2})_0=R\Ga(E,\GG)
\end{equation}
where $E=\pi^{-1}(0)$ is the exceptional divisor,
$$\GG=\wt{j}_*(L(\psi(\frac{1}{y^n})\frac{\chi_1}{\chi_2}(x)
\chi_2(y)))|_E.$$
Consider the open chart in $\wt{\A^2}$ with coordinates $x, u$ such that
$y=xu$. In this chart $E$ is given by the equation $x=0$ while
$u$ is the local coordinate on $E$. We have
$$\psi(\frac{1}{y^n})\frac{\chi_1}{\chi_2}(x)\chi_2(y)=
\psi(\frac{1}{x^nu^{n}})\chi_1(x)\chi_2(u).$$
Thus, for $u\neq 0$ we have $\GG_u=0$ while 
$\GG_0=(\wt{F}^{-n,-n}_{\chi_1,\chi_2})_0$. Now (\ref{rga}) implies that
$\GG_0$ is a direct summand in     
$(\wt{F}^{0,-n}_{\frac{\chi_1}{\chi_2},\chi_2})_0=0$. Therefore,
$\GG_0=(\wt{F}^{-n,-n}_{\chi_1,\chi_2})_0=0$. 

\noindent
(ii) Consider the standard action of $\G_m$ on $\A^2$. Then with respect to
an action of $t\in\G_m$ the sheaf $\wt{F}^{n,n}_{\chi_1,\chi_2}$
gets tensored with $(L_{\chi_1\chi_2^{-1}})_t$. Since $\chi_1\chi_2^{-1}$
is non-trivial this immediately implies the result.

\noindent
(iii) Consider first the case $m=n=d$. 
Then using the above notation we can write 
$$(\wt{F}^{d,d}_{\chi,\chi^{-1}})_0=R\Ga(E,\wt{j}_*(L(\psi(\frac{x^d}{y^d})
\chi(\frac{x}{y})))|_E).$$
On the open chart $U_1\sub\wt{\A^2}$ with coordinates $x, u$ such that
$y=xu$ we have 
$$L(\psi(\frac{x^d}{y^d})\chi(\frac{x}{y}))=L(\psi(\frac{1}{u^d})
\chi^{-1}(u)).$$
On the second open chart $U_2\sub\wt{\A^2}$
with coordinates $v, y$ such that $x=vy$ we have
$$L(\psi(\frac{x^d}{y^d})\chi(\frac{x}{y}))=L(\psi(v^d)
\chi(v)).$$
Hence, 
$$(\wt{F}^{d,d}_{\chi,\chi^{-1}})_0=R\Ga(E,L(\psi(v^d)\chi(v)))$$
where the sheaf $L(\psi(v^d)\chi(v))$ is extended by zero from
$\A^1$ to $\P^1=E$. This implies immediately that
$$(F^{d,d}_{\chi,\chi^{-1}})_0=H^1_c(\A^1,L(\psi(v^d)\chi(v)))[1].$$

Now in the general case we have
$$(F^{m,-n}_{\chi^{m'},\chi^{-n'}})_0=
H^0(G,(F^{m'n'd,m'n'd}_{\chi^{m'n'},\chi^{-m'n'}})_0)=
H^0(G,H^1_c(\A^1,L(\psi(v^{m'n'd})\chi(v^{m'n'}))))[1]$$
The latter space up to a shift of degree has form 
$H^0(G, R\Ga_c(f^*L(\psi(v^d))\chi(v)))$ where
$f:\A^1\ra\A^1:v\mapsto v^{m'n'}$ is a covering with Galois group $G$.
Therefore, this space is isomorphic to $R\Ga_c(L(\psi(v^d)\chi(v)))$.
\ed

The theorem \ref{2dim} admits the following partial generalization to higher
dimensions.

\begin{thm}\label{ndim} 
Let $n_i>0$, $i=1,\ldots,k$, $m_j$, $j=1,\ldots,l$ and $n>0$
be integers which are prime to $q$.
Set $m_j=m_j'd$ where $d=gcd(m_1,\ldots,m_l)$.
Then

\noindent
(i) $(\wt{F}^{-n_1,\ldots,-n_k}_{\chi_1,\ldots,\chi_k})_0=
(F^{-n_1,\ldots,-n_k}_{\chi_1,\ldots,\chi_k})_0=0$
for any $\chi_1,\ldots,\chi_k$.

\noindent
(ii) $(\wt{F}^{m_1,\ldots,m_l}_{\eta_1,\ldots,\eta_l})_0=
(F^{m_1,\ldots,m_l}_{\eta_1,\ldots,\eta_l})_0=0$
unless there exists a character $\eta$ such that $\eta_i=\eta^{m'_i}$
for all $i$.

\noindent
(iii) $(\wt{F}^{n,-n_1,\ldots,-n_k}_{1,1,\ldots,1})_0=0$
provided that $gcd(n,n_1,\ldots,n_k)=1$.

\end{thm}

\Pf . (i) Considering the covering of $\A^{k+1}$ induced by
$(x_1,x_2)\mapsto (x_1^{n_2},x_2^{n_1})$ we reduce to the case $n_1=n_2$.
Then the argument with blow-up along $x_1=x_2=0$ similar to the case (i)
of Theorem (\ref{2dim}) allows the induction in $k$.

\noindent
(ii) Consider the natural action of the torus
$$T=\{(t_1,\ldots,t_l)\in\G_m^l:\prod_{i=1}^l t_i^{m'_i}=1\}$$
on $\A^l$. Then under the action of $(t_1,\ldots,t_l)\in T$
both the sheaves $\wt{F}^{m_1,\ldots,m_l}_{\eta_1,\ldots,\eta_l}$
and $F^{m_1,\ldots,m_l}_{\eta_1,\ldots,\eta_l}$ get tensored
with the space $(L_{\eta_1})_{t_1}\otimes\ldots\otimes(L_{\eta_l})_{t_l}$.
Now unless there exists $\eta$ such that
$\eta_j=\eta^{m'_j}$ for all $j$
we can find a one-parameter subgroup $i:\G_m\ra T$ such that
$i^*(L_{\eta_1}\otimes\ldots L_{\eta_l})\simeq L_{\chi}$ where
$\chi$ is a non-trivial character. Then under the action of
$i(t)$ our sheaves get tensored with $(L_{\chi})_t$ 
which immediately implies the vanishing of their stalks at zero.

\noindent
(iii) We proceed by induction in $k$. The case $k=1$ follows from Theorem
\ref{2dim}, (iii). Now assume that $k>1$.
Let us denote the variables by $(x,x_1,\ldots,x_k)$. Let
$n=n'd$, $n_1=n'_1d$ where $d=gcd(n,n_1)$.
Consider first the Galois covering $(\ov{x},\ov{x_1})\mapsto
(x=\ov{x}^{n'_1},x_1=\ov{x_1}^{n'})$ with Galois group
$G=\mmu_{n'n'_1}$. Then we have
$$(\wt{F}^{n,-n_1,\ldots,-n_k}_{1,1,\ldots,1})_0=H^0(G,
(Rj_*L(\psi(\frac{x^{dn'n'_1}}{x_1^{dn'n'_1}x_2^{n_2}\ldots})))_0).$$
Let
$\pi:\wt{\A^{k+1}}\ra \A^{k+1}$ be the blow-up of $\A^{k+1}$
along the subspace
$x=x_1=0$. Let $\wt{j}:\A^1\times\G_m^k\ra \wt{\A^{k+1}}$
be the natural embedding where the first factor $\A^1$ corresponds
to the variable $x$. We have
$$(Rj_*L_{\psi}(\frac{x^{dn'n'_1}}{x_1^{dn'n'_1}x_2^{n_2}\ldots}))_0
=R\Ga(E_0,
R\wt{j}_*L_{\psi}(\frac{x^{dn'n'_1}}{x_1^{dn'n'_1}x_2^{n_2}\ldots})|_{E_0})$$
where $E_0=\pi^{-1}(0)\simeq\P^1$. 

On the open chart $U_1\sub\wt{\A^{k+1}}$ with coordinates 
$x, u, x_2,\ldots$
such that $x_1=xu$ we have 
$$L_{\psi}(\frac{x^{dn'n'_1}}{x_1^{dn'n'_1}x_2^{n_2}\ldots x_k^{n_k}})=
L_{\psi}(\frac{1}{u^{dn'n'_1}x_2^{n_2}\ldots x_k^{n_k}})$$
The extension of this sheaf has zero stalks over $E_0\cap U_1$
by (i). On the second open chart $U_2\sub\wt{\A^{k+1}}$
with coordinates $v,x_1,x_2,\ldots$ such that $x=vx_1$ we have
$$L_{\psi}(\frac{x^{dn'n'_1}}{x_1^{dn'n'_1}x_2^{n_2}\ldots x_k^{n_k}})=
L_{\psi}(\frac{v^{dn'n'_1}}{x_2^{n_2}\ldots x_k^{n_k}}).$$
Thus, the sheaf
$$R\wt{j}_*L_{\psi}(\frac{x^{dn'n'_1}}{x_1^{dn'n'_1}x_2^{n_2}\ldots 
x_k^{n_k}})|_{E_0}$$
is supported at one point $v=0$ and we get
$$(Rj_*L_{\psi}(\frac{x^{dn'n'_1}}{x_1^{dn'n'_1}x_2^{n_2}\ldots x_k^{n_k}}))_0=
(Rj_*L_{\psi}(\frac{v^{dn'n'_1}}{x_2^{n_2}\ldots x_k^{n_k}}))_0.$$
Hence,
$$(\wt{F}^{n,-n_1,\ldots,-n_k}_{1,1,\ldots,1})_0=H^0(G,
(Rj_*L_{\psi}(\frac{v^{dn'n'_1}}{x_2^{n_2}\ldots x_k^{n_k}}))_0)=
(Rj_*L_{\psi}(\frac{v^d}{x_2^{n^2}\ldots x_k^{n_k}}))_0.$$
Since $gcd(d,n_2,\ldots,n_k)=gcd(n,n_1,\ldots,n_k)=1$ this space is zero
by induction assumption.
\ed

\subsection{Traces of Frobenius on the stalks of
Goreski-MacPherson extensions}\label{trfr}

In this section we show how to compute the traces of $\Frob_q$
on the stalks of the sheaves $F^{n_1,\ldots,n_k}_{\chi_1,\ldots,\chi_k}$
on $\A^k$ introduced above. 
It suffices to consider the stalk at zero, and according
to Theorem \ref{ndim},(ii) the non-zero stalk at zero can appear only
for the sheaf of the form 
$F^{dn_1,\ldots,dn_k}_{\chi^{n_1},\ldots,\chi^{n_k}}$.
First let us consider the particular case when
$n_i=\pm 1$ for all $i$.
In order to formulate the answer we have to introduce two families
of polynomials in $q$:
\begin{equation}\label{anm}
a(n,m)=\cases \sum_{i=0}^{m-1}{m-1\choose i}{n-1\choose i+1}q^{i+1}, & m\ge 1
              \\              1, & m=0\endcases 
\end{equation}
\begin{equation}\label{bnm}
b(n,m)=\sum_{i=0}^{m-1}{m-1\choose i}{n-1\choose i}q^i.
\end{equation}
These polynomials satisfy the following recursive relations
$$b(n,m)-b(n,m-1)=a(n,m-1),$$
$$a(n,m)-a(n-1,m)=qb(n-1,m).$$
For a local system $L$ on $U\stackrel{j}{\hra}\A^N$ let us denote
$$j_{!*}L:=j_{!*}(L[N])[-N].$$

\begin{thm}\label{GMtr} 
Let $\chi$ be a non-trivial character, $d\ge 1$. Then
\begin{equation}\label{trfr1}
\Tr(\Frob_q,j_{!*}L(\psi(\frac{x_1^d\ldots x_n^d}{y_1^d\ldots y_m^d}))_0)
=a(n,m)+b(n,m)\cdot\sum_{t\in\F_q}\psi(t^d),
\end{equation}
\begin{equation}\label{trfr2}
\Tr(\Frob_q,j_{!*}L(\psi(\frac{x_1^d\ldots x_n^d}{y_1^d\ldots y_m^d})\chi
(\frac{x_1\ldots x_n}{y_1\ldots y_m}))_0)
=b(n,m)\cdot\sum_{t\in\F_q^*}\psi(t^d)\chi(t).
\end{equation}
\end{thm}

\Pf . Consider the blow-up $\pi:\wt{\A^{n+m}}\ra\A^{n+m}$ along
the subspace $x_1=y_1=0$. Let $\wt{j}$ be the natural embedding of
$\G_m^{n+m}$ into $\wt{\A^{n+m}}$. We claim that
\begin{equation}\label{blowup}
R\pi_*(\wt{j}_{!*}L(\psi(\frac{x_1^d\ldots x_n^d}{y_1^d\ldots y_m^d})))
\end{equation}
is the Goreski-MacPherson extension from $\G_m^{n+m}$, and the similar
statement holds for 
$L(\psi(\frac{x_1^d\ldots x_n^d}{y_1^d\ldots y_m^d})
\chi(\frac{x_1\ldots x_n}{y_1\ldots y_m}))$. 
Let us denote by $S(k)$ the subset of $\A^{n+m}$ where exactly
$k$ coordinates vanish (so $k\ge 2$).
It suffices to check that the stalks of the sheaf (\ref{blowup}) over
points of $S(k)$ with $x_1=y_1=0$ 
are concentrated in degrees $<k$. Now for 
$p=(0,a_2,\ldots,a_n,0,b_2,\ldots,b_m)\in S(k)$ we have          
$$
R\pi_*(\wt{j}_{!*}L(\psi(\frac{x_1^d\ldots x_n^d}{y_1^d\ldots y_m^d})))_p=
R\Ga(\pi^{-1}(p), \wt{j}_{!*}
L(\psi(\frac{x_1^d\ldots x_n^d}{y_1^d\ldots y_m^d}))|_{\pi^{-1}(p)}).
$$
We can cover $\wt{\A^{n+m}}$ by two open affine charts:
$U_1$ with coordinates $(u,x_2,\ldots,x_n,y_1,\ldots,y_m)$ such that
$x_1=uy_1$ and $U_2$ with coordinates
$(x_1,\ldots,x_n,v,y_2,\ldots,y_m)$ such that $y_1=vx_1$. 
Both $U_1$ and $U_2$ are isomorphic to $\A^{n+m}$ and we will denote
by $S(k)$ the strata given by vanishing of coordinates in $U_1$ and 
$U_2$. Let us represent $\pi^{-1}(p)$ as the disjoint union of 
$\pi^{-1}(a)\cap U_1\cap U_2\simeq\G_m$ and two points $p_1$, $p_2$,
where $p_1\in U_1$ has coordinates $(u=0,a_2,\ldots,a_n,0,b_2,\ldots,b_m)$
and $p_2\in U_2$ has coordinates $(0,a_2,\ldots,a_n,v=0,b_2,\ldots,b_m)$.
We have $\pi^{-1}(p)\cap U_1\cap U_2\subset S(k-2)$ 
while $p_1$ and $p_2$ are in $S(k-1)$. Consider the exact triangle
\begin{align*}
&R\Ga(\pi^{-1}(p)\cap U_1\cap U_2, \wt{j}_{!*}
L(\psi(u^d\frac{x_2^d\ldots x_n^d}{y_2^d\ldots y_m^d})))\ra
R\Ga(\pi^{-1}(p), \wt{j}_{!*}
L(\psi(\frac{x_1^d\ldots x_n^d}{y_1^d\ldots y_m^d})))\\
&\ra
\oplus_{i=1,2} 
(\wt{j}_{!*}L(\psi(\frac{x_1^d\ldots x_n^d}{y_1^d\ldots y_m^d})))_{p_i}\ra
\ldots
\end{align*}
Assume first that $k>2$. Then making the
change of variables $x_2'=ux_2$ (and leaving all the other variables the same)
we can rewrite the first term of the above triangle as
$$R\Ga_c(\G_m,\ov{\Q}_l)\otimes
j_{!*}L(\psi(\frac{(x'_2)^d\ldots x_n^d}{y_2^d\ldots y_m^d}))_{p'}$$
where $p'=(a_2,\ldots,a_n,b_2,\ldots,b_m)$. Since $p'\in S(k-2)$
and $R\Ga_c(\G_m,\ov{\Q}_l)$ lives in degrees $1$ and $2$, this
term is concentrated in degrees $<k$. Now assume that $k=2$, i.e.
$a_2\ldots a_nb_2\ldots b_m\neq 0$. Then we can make the change
of variables $u'=ux_2\ldots x_ny_2^{-1}\ldots y_m^{-1}$ so that the first
term of the above triangle takes form
$$R\Ga(\G_m,L_{\psi}(u^d))$$
which is concentrated in degree $1<k=2$. On the other hand, since
$p_i\in S(k-1)$ the last
term of the above triangle is concentrated in degrees $<k-1$. Therefore,
our claim follows, so the sheaf (\ref{blowup}) is the Goreski-MacPherson
extension. The similar argument works for the sheaf involving the non-trivial
character $\chi$.

Let us denote
$$A(n,m):=
\Tr(\Frob_q,j_{!*}L(\psi(\frac{x_1^d\ldots x_n^d}{y_1^d\ldots y_m^d}))_0).$$
Then for $n+m>2$ the above exact triangle shows
that
\begin{equation}\label{rec}
A(n,m)=(q-1)A(n-1,m-1)+A(n-1,m)+A(n,m-1).
\end{equation}
Notice that $a(n,m)$ and $b(n,m)$ are solutions of (\ref{rec}) with
initial values $a(0,m)=b(0,m)=0$, $a(n,0)=1$, $a(1,1)=0$,
$b(n,0)=0$, $b(1,1)=1$. Since 
$A(1,1)=\sum_{t\in\F_q}\psi(t^d)$ by Theorem \ref{2dim},(iii),
$A(0,m)=0$ by Theorem \ref{ndim},(i) and $A(n,0)=1$ this proves
(\ref{trfr1}). The proof of (\ref{trfr1}) is similar.
\ed 

In the general case we can proceed as follows. The stalk
at zero of the sheaf
\begin{equation}\label{sheaf}
j_{!*}L(\psi(\frac{x_1^{n_1d}\ldots x_k^{n_kd}}{y_1^{m_1d}\ldots
y_l^{m_ld}})\chi(\frac{x_1^{n_1}\ldots x_k^{n_k}}{y_1^{m_1}\ldots
y_l^{m_l}}))
\end{equation}
where $n_i>0$, $m_j>0$, $gcd(n_i,q)=gcd(m_j,q)=1$, can be computed
as $G$-invariants of the stalk at zero of the sheaf
$$j_{!*}L(\psi(\frac{x_1^{Nd}\ldots x_k^{Nd}}{y_1^{Nd}\ldots
y_l^{Nd}})\chi(\frac{x_1^N\ldots x_k^N}{y_1^N\ldots
y_l^N}))$$
where $N=lcm(n_1,\ldots,n_k,m_1,\ldots,m_l)$ (by $lcm$ we mean the
least common multiple), $G$ is the product of the groups
$\mmu(N/n_1)$, ..., $\mmu(N/m_l)$ acting on the variables in the natural
way. Now we can use the blow-up along $x_1=y_1=0$ and 
the $G$-equivariant exact triangle similar
to the one considered in the above proof. As a result we get the following
recursion relation for the trace of $\Frob_q$ acting on the stalk of
(\ref{sheaf}) at zero:
$$a^{n_1,n_2,\ldots,n_k}_{m_1,m_2,\ldots,m_l}=
(q-1)a^{n'_2,\ldots,n_k}_{m_2,\ldots,m_l}+a^{n,n_2,\ldots,n_k}_{m_2,\ldots,m_l}
+a^{n_2,\ldots,n_k}_{n,m_2,\ldots,m_l}$$
for $k+l>2$ where
$n'_2=lcm(n_1,m_1,n_2)$, $n=lcm(n_1,m_1)$.
Using this recursion relation and Theorem \ref{2dim} we can in principle
compute all these traces.

\subsection{Pointwise purity}

Recall (see \cite{DeWeil}) that an object $K\in\D^b_c(\A^n)$ is called
{\it pointwise pure of weight} $w$
if for every closed point $x\in\A^n$ such that
$k(x)=\F_{q_1}$ the endomorphism $\Frob_{q_1}$ acts on $H^iK|_x$
with eigenvalues which are algebraic numbers with all conjugates
of absolute value $q_1^{\frac{w+i}{2}}$.

\begin{lem} Let $K$ be a $\G_m$-equivariant sheaf on $\A^n$ where $\G_m$
acts on $\A^n$ by
$$t(x_1,\ldots,x_n)=(t^{d_1}x_1,\ldots,t^{d_n}x_n).$$
Assume that all weights $d_i$ are positive integers. Then
the natural map
$$R\Ga(\A^n,K)\ra K|_0$$
is an isomorphism.
\end{lem}

\Pf . Consider the coordinate stratification of $\A^n$. It suffices
to prove that for any stratum $S\subset\A^n\setminus\{0\}$
and a $\G_m$-equivariant sheaf $K$ on $S$ one has
$H^q(\A^n,j_!K)=0$ for any $q$ where $j:S\ra\A^n$ is the embedding.
Without loss of generality we can assume that $S$ is the open stratum:
$S=\G_m^n$. Let $d$ be the greatest common divisor of $d_1,\ldots,d_n$.
Consider the covering 
$$\pi:\A^n\ra\A^n:(x_1,\ldots,x_n)\mapsto (x_1^d,\ldots,x_n^d).$$
Set $K'=(\pi|_{\G_m^n})^*K$. Then $K'$ is a $\G_m$-equivariant sheaf
on $\G_m^n$ with respect to the action with the weights $(d_1/d,\ldots,
d_n/d)$. Since $j_!K$ is the direct summand in $\pi_*(j_!K')$ it suffices
to prove that cohomologies of $j_!K'$ vanish. Thus, we can assume from
the beginning that $d=1$. Then the action of $\G_m$ on $\G_m^n$ is free.
Let $p:\G_m^n\ra T$ be the quotient under this action (so that $T$ is a torus).
We have $K=p^*L$ for some sheaf $L$ on $T$. Let
$f:\wt{\A^n}\ra \A^n$ be the weighted blow-up of $\A^n$ along the origin,
i.e. $\wt{\A^n}=\Proj k[x_1,\ldots,x_n,x_1t^{d_1},\ldots,x_nt^{d_n}]$
where $x_i$ are coordinates on $\A^n$ ($\deg x_i=0$), $t$ is an
independent variable of degree $1$.    
The morphism $f$ is proper and the inclusion $j:\G_m^n\ra\A^n$ factors through
the inclusion $\wt{j}:\G_m^n\ra\wt{\A^n}$. Hence, it suffices to prove
that $H^q(\wt{\A^n},\wt{j}_!K)=0$. 
Note that there is a natural
projection $\wt{p}:\wt{\A^n}\ra\P(d_1,\ldots,d_n)$
where $\P(d_1,\ldots,d_n)=\Proj k[y_1,\ldots,y_n]$ ($\deg y_i=d_i$) is the 
corresponding weighted projective space.
We can identify $T$ with an open subset of $\P(d_1,\ldots,d_n)$ defined by
$y_1\ldots y_n\neq 0$ so that the following diagram is commutative.
\begin{equation}
\begin{array}{ccccc}
\G_m^n &\lrar{} & \wt{p}^{-1}(T) &\lrar{} &\wt{\A^n}\\
&\ldrar{p} &\ldar{} &&\ldar{\wt{p}}\\
&& T &\lrar{} & \P(d_1,\ldots,d_n)
\end{array}
\end{equation}
Choosing a section of the homomorphism of tori $p:\G_m^n\ra T$ we can lift
the natural action of $T$ on $\P(d_1,\ldots,d_n)$ to an action of $T$ on
$\wt{\A^n}$. Hence, the projection $\wt{p}^{-1}(T)\ra T$ is
a locally trivial fibration with fiber $F$ which is equal to a generic 
$\G_m$-orbit on $\A^n$: $F=\{(t^{d_1},\ldots,t^{d_n}),t\in\A^1\}$.
Furthermore there is a canonical zero
section $\sigma:T\ra\wt{p}^{-1}(T)$ and 
$\G_m^n\subset\wt{p}^{-1}(T)$ is the complement to $\sigma(T)$.
Now the fact that $H^q(F,k_!\ov{\Q}_{l,F-0})=0$ implies easily
that $R^q\wt{p}_*(\wt{j}_!K)=0$.
\ed

\begin{prop}\label{pure}
The perverse sheaves $F^{n_1,\ldots,n_k}_{\chi_1,\ldots,\chi_k}$
are pointwise pure.
\end{prop} 

\Pf . It suffices to prove that the stalk at zero of the sheaf (\ref{sheaf})
is pure of weight $0$.
By Galois covering argument it suffices to prove the purity of the stalk
at zero of the sheaf 
$$K=j_{!*}L(\psi(\frac{x_1^d\ldots x_n^d}{y_1^d\ldots y_m^d})\chi
(\frac{x_1\ldots x_n}{y_1\ldots y_m}))_0$$
where $m,n>0$.
Note that $K$ is pure of weight $0$ (see \cite{BBD}, 5.3.2). In particular,  
$K|_0$ is of weight $\le 0$.
On the other hand, by the principal theorem of \cite{DeWeil} (\cite{DeWeil},
3.3.1, 6.2.3) the complex $R\Ga(\A^n,K)$ is of weight $\ge 0$.
Applying the above lemma to the $\G_m$-action
$$t(x_1,\ldots,x_n,y_1,\ldots,y_m)=(t^mx_1,\ldots,t^mx_n,t^ny_1,\ldots,
t^ny_m)$$
we deduce that $K|_0$ is pure of weight $0$.
\ed

\begin{cor} Let $l$ be the number of negative integers among
$(n_1,\ldots,n_k)$, $N$ be the least common multiple of 
$(|n_1|,\ldots,|n_k|)$. 
Then the eigenvalues of $\Frob_q$ on the stalk of the sheaf
$$j_{!*}L(\psi(x_1^{dn_1}\ldots x_k^{dn_k})\chi(x_1^{n_1}\ldots x_k^{n_k}))$$
at $0$ for $\chi\neq 1$ have form $q^i\la$ 
where $i\in\Z$, $0\le i< N$, and 
$\la^r=-g(\eta)$ for some character $\eta$ of $\F_{q^r}^*$
satisfying $\eta^{dN}=\chi^N\circ\Nm_r$.
If $\chi=1$ in addition the eigenvalues of the form $q^{i+1}$ for $0\le i<N$
can appear.
\end{cor}

\section{Identities between Gauss sums}
\label{gauss-sec}

\subsection{Relations between cyclothomic divisors}\label{cycl}

For every $N\in\Z_{>0}$ let us denote by $A_N$ the abelian group generated by 
symbols $[s,n]_N$ where $s\in\Z/N\Z$, $n\in\Z_{>0}$.
subject to the following relations: for every $d|N$
we have
\begin{equation}\label{defrel}
[ds,dn]_N=\sum_{i=0}^{d-1}[s+i\frac{N}{d},n]_N.
\end{equation}
Also for every prime $p$ we can consider the group $A_N^{(p)}$ defined
in the same way as $A_N$ except that we allow only symbols $[s,n]_N$
with $(n,p)=1$, and the relation (\ref{defrel}) is imposed for every
$d|N$ such that $(d,p)=1$.

\begin{lem}\label{basis}
The elements $[s,n]_N$ such that $gcd(s,n,N)=1$ form a basis of
$A_N$ (resp. $A_N^{(p)}$). In particular, $A_N^{(p)}$ is a subgroup
in $A_N$.
\end{lem} 

The proof is straightforward and is left to the reader.

For every set $S$ let us denote by $\Div(S)$ the group $\oplus_S \Z$,
i.e. the group of formal linear combinations of elements of $S$
with integer coefficients. We'll call elements of $\Div(S)$ divisors on $S$. 
In particular, we want to consider
$\Div(\Q/\Z)$. For every pair $(r,n)$ where $r\in\Q$, $n\in\Z_{>0}$,
we consider the divisor 
\begin{equation}\label{drn}
D_{r,n}=(r)+(r+1/n)+\ldots (r+(n-1)/n)
\end{equation}
on  
$\Q/\Z$. The divisor $D_{r,n}$ depends only on $n$ and on the residue
class of $r$ modulo $\frac{1}{n}\Z$.
We have the following relations between these divisors:
$$D_{r,dn}=\sum_{i=0}^{d-1} D_{r+i/dn,n}.$$
This means that we can define the homomorphism 
$$\a_N:A_N\ra\Div(\Q/\Z)$$
by the formula $\a_N([s,n]_N)=D_{\frac{s}{nN},n}$.

One immediately checks that for every $M\in\Z_{>0}$ there is a homomorphism 
$\phi_{M,N}:A_N\ra A_{MN}$ sending $[s,n]_N$ to $[Ms,n]_{MN}$.
Moreover, one has $\phi_{K,MN}\circ\phi_{M,N}=\phi_{KM,N}$ and
$\a_{N}=\a_{MN}\circ\phi_{M,N}$.

\begin{lem}\label{dirsum}
The homomorphism $\phi_{M,N}$ identifies $A_N$ (resp. $A^{(p)}_N$) with a
direct summand of $A_{MN}$ (resp. $A^{(p)}_{MN}$).
\end{lem}

\Pf . It suffices to consider the case when $M$ is prime.
Let us look at the images of basis elements $[s,n]_N$,
$gcd(s,n,N)=1$. If $gcd(n,M)=1$ then
$\phi_{M,N}([s,n]_N)=[Ms,n]_{MN}$ is a basis element in $A_{MN}$.
Otherwise, $M|n$ and we have
$$\phi_{M,N}([s,n]_{N})=\sum_{i=0}^{M-1}[s+iN,\frac{n}{M}]_{MN}.$$
Since $gcd(s,n,N)=1$ it follows that $gcd(s,M,N)=1$, therefore,
replacing $s$ by $s+N$ if necessary we can assume that $gcd(s,M)=1$.
Then $[s,\frac{n}{M}]_{MN}$ is a basis element in $A_{MN}$.
Clearly, the basis elements obtained in this way are all different
which implies our statement.
\ed

\begin{thm}\label{ker} One has $\ker(\a_N)=0$.
\end{thm}
 
Assume that $\a_N(x)=0$ where
$x=\sum_i m_i [s_i,n_i]_N$ with $gcd(s_i,n_i,N)=1$.
Set $M=\prod_i n_i$. We claim that
$\phi_{M,N}(x)=0$. Indeed,
we have
\begin{equation}\label{phiMN}
\phi_{M,N}([s_i,n_i]_N)=[Ms_i,n_i]_{MN}=\sum_{j=0}^{n_i-1}
[M\wt{s}_i+j\frac{MN}{n_i},1]_{MN}
\end{equation}
where $\wt{s}_i\in\Z$ is a representative of $s_i$.
By Lemma \ref{basis} the map $\a_{MN}$ is injective on the subgroup
generated by elements $[t,1]_{MN}$. Therefore, $\phi_{M,N}(x)=0$.
Hence, $x=0$ by Lemma \ref{dirsum}.
\ed

\subsection{Application to Gauss sums}

One can generalize the content of \ref{cycl} to the case
of cyclic groups without fixed generators. Namely, for every
finite cyclic group $G$ we can define the abelian group
$A(G)$ generated by symbols $[g,n]$ where $g\in G$, $n\in\Z_{>0}$
subject to relations
\begin{equation}\label{defrel2}
[g^d,dn]=\sum_{h\in H}[gh,n]
\end{equation}
for every subgroup $H\sub G$, where $d=|H|$.
Given a prime number $p$ we can define similarly the group
$A^{(p)}(G)$ using only the symbols $[g,n]$ with $gcd(n,p)=1$. 
Thus, we have $A_N=A(\Z/N\Z)$, $A_N^{(p)}=A^{(p)}(\Z/N\Z)$.

Note that Lemma \ref{dirsum} can be reformulated as follows:
for every inclusion of finite cyclic groups $H\subset G$
the induced homomorphism $A(H)\ra A(G)$ identifies $A(H)$ with
the direct summand of $A(G)$. Similar property holds for the homomorphism
$A^{(p)}(H)\ra A^{(p)}(G)$.

Now let us fix a finite field $\F_q$, where $q=p^s$.
For every $d>0$ we denote
by $X_d=X(\F_{q^d}^*)$ the group of characters of $\F_{q^d}^*$. For every
$d_1|d_2$ we have the inclusion $X_{d_1}\ra X_{d_2}$ induced
by the norm homomorphism $\F_{q^{d_2}}^*\ra\F_{q^{d_1}}^*$.
Let us denote by $X$ the direct limit of the system $(X_d)$ with respect to
these inclusions. The group $X$ is isomorphic non-canonically
to the $q$-prime part of $\Q/\Z$.

Let us define the homomorphism $\a^d:A^{(p)}(X_d)\ra\Div(X)$
by sending $[\chi,n]$, $\chi\in X_d$, $n\in\Z_{>0}$, to the divisor
\begin{equation}\label{Dchin}
D_{\chi,n}=\sum_{\xi\in X:\xi^n=\chi}(\xi).
\end{equation}
Then Theorem \ref{ker} implies that all the homomorphisms $\a^d$
are injective. Moreover, the induced injective homomorphism
\begin{equation}\label{isdiv}
\lim_d A^{(p)}(X_d)\ra\Div(X)
\end{equation} 
is clearly surjective, therefore, it is an
isomorphism.

Let us fix a non-trivial additive character $\psi:\F_q\ra\C^*$.
Then to every generator $[\chi,n]$ of $A^{(p)}(X(\F_q^*))$,
we can associate the following function $f_{\chi,n}$
on $X(\F_q^*)$:
$$f_{\chi,n}(\la)=\frac{g(\la^n\chi)}{\la(n^n)g(\chi)}$$
(where by our convention $g(1)=-1$).
It is easy to see that Hasse-Davenport formula (\ref{HasDav}) implies
that the map $[\chi,n]\mapsto f_{\chi,n}$ extends to a homomorphism
$A^{(p)}(X(\F_q^*))\ra\CC(X(\F_q^*),\C^*)$ where 
$\CC(S,\C^*)=(\C^*)^S$ is the group of $\C^*$-valued functions on $S$.

It is convenient to set 
\begin{equation}\label{Dchin2}
D_{\chi^{-1},-n}=-D_{\chi,n}. 
\end{equation}
We have the following
corollary of Theorem \ref{ker}.

\begin{cor}\label{maincor} 
Assume that $\sum_{i=1}^k D_{\chi_i,n_i}=0$
in $\Div(X)$, where $\chi_i\in X(\F_q^*)$, $gcd(n_i,q)=1$.
Then for every character $\la$ of $\F_q^*$ one has
$$\prod_i\frac{g(\la^{n_i}\chi_i)}{\la(n_i^{n_i})g(\chi_i)}=
q^{m(\la)}$$
for some $m(\la)\in\Z$.
In particular, if $\la^{n_i}\chi_i\neq 1$ for all $i$ then
$2m(\la)$ is the number of $i$ such that $\chi_i=1$.
\end{cor}

\Pf . This follows from the fact that $g(\la^{-1})\equiv\la(-1)g(\la)^{-1}
\mod q^{\Z}$ for any $\la$.
\ed

In fact, the above corollary in some sense describes all multiplicative
identities between Gauss sums which hold universally over all extensions
of a given finite field. Here is a more precise statement.

\begin{thm} Let $(\chi_1,\ldots,\chi_k)$ be a collection of characters of
$\F_q^*$, $(n_1,\ldots,n_k)$ be integers such that $gcd(n_i,q)=1$.
Assume that for any $d\ge 1$ and for any character
$\la\in X(\F_{q^d}^*)$ such that $\la^{-1}\neq \chi_i\circ\Nm_d$ for any $i$, one
has
$$
\prod_{i=1}^k (-g(\la^{n_i}(\chi_i\circ\Nm_d),\psi\circ\Tr_d))=
c^d\cdot\la(a)
$$
for some constants $c\in\C^*$, $a\in\F_q^*$. Then
$\sum_{i=1}^k D_{\chi_i,n_i}=0$.
\end{thm}

\Pf . Since $D_{\chi_i\circ\Nm_d,n_i}=D_{\chi,n_i}$ we can pass to any
extension $\F_{q^d}$ of $\F_q$. Thus, we can assume that $\chi_i$ is
$n_i$-th power of some character and applying Hasse-Davenport identity
(\ref{HasDav}) we reduce ourselves to the case $n_i=\pm 1$.
In other words, it suffices to prove that if for two
collections of characters of $\F_q^*$: $(\chi_1,\ldots,\chi_k)$
and $(\eta_1,\ldots,\eta_l)$ and for some constants $c\in\C^*$ and
$a\in\F_q^*$ one has
$$
\prod_{i=1}^k (-g(\la(\chi_i\circ\Nm_d),\psi\circ\Tr_d))=c^d\cdot\la(a)\cdot
\prod_{j=1}^l (-g(\la(\eta_j\circ\Nm_d),\psi\circ\Tr_d))
$$
for all $\la\in X(\F_{q^d}^*)$ such that $\la^{-1}$ is different
from all $\chi_i\circ\Nm_d$ and $\eta_j\circ\Nm_d$,
then $\chi_i=\eta_j$ for some $(i,j)$. Let
$$\{\chi_1,\ldots,\chi_k,\eta_1,\ldots,\eta_l\}=\{\mu_1,\ldots\mu_r\}$$
where the characters $\mu_i$ are all different.
Then we can rewrite the above identity as follows
\begin{align*}
&\prod_{i=1}^k (-g(\la(\chi_i\circ\Nm_d),\psi\circ\Tr_d))-\sum_{i=1}^rc_i^d
\de(\la(\mu_i\circ\Nm_d))=\\
&c^d\cdot\la(a)\cdot
\left(\prod_{j=1}^l (-g(\la(\eta_j\circ\Nm_d),\psi\circ\Tr_d))-
\sum_{i=1}^rb_i^d\de(\la(\mu_i\circ\Nm_d))\right)
\end{align*}
for some constants $c_i$, $b_i$, where $\de$ is the delta-function
at the trivial character.
Let us denote 
$$K_{\psi;\chi_1,\ldots,\chi_k}(t)=\sum_{x_1\ldots x_k=t}\psi(x_1+\ldots+x_k)
\chi_1(x_1)\ldots\chi_k(x_k)$$
for $t\in\F_q^*$. Then the multiplicative Fourier transform of 
$K_{\psi;\chi_1,\ldots,\chi_k}$ is the function
$$\la\mapsto \prod_{i=1}^k g(\la\chi_i).$$
Thus, applying the inverse Fourier transform to the above equality
we obtain the following equality of functions on $\F_{q^d}^*$:
\begin{align*}
& (-1)^kK_{\psi\circ\Tr_d;\chi_1\circ\Nm_d,\ldots,\chi_k\circ\Nm_d}-
\frac{1}{q-1}\sum_{i=1}^r c_i^d\mu_i\circ\Nm_d=\\
&c^d t_{a^{-1}}^*
\left((-1)^lK_{\psi\circ\Tr_d;\eta_1\circ\Nm_d,\ldots,\eta_l\circ\Nm_d}-
\frac{1}{q-1}\sum_{i=1}^r b_i^d\mu_i\circ\Nm_d\right)
\end{align*}
where $t^*_{a^{-1}}f(x):=f(a^{-1}x)$.
According to Theorem 7.8 of \cite{De1} there exists an irreducible
local system on $\G_m$ whose trace functions over extensions of $\F_q$
are given by $K_{\psi\circ\Tr_d;\chi_1\circ\Nm_d,\ldots,\chi_n\circ\Nm_d}$.
Therefore, by Theorem 1.1.2 of \cite{La2}
the above equality implies the similar equality in the
Grothendieck group of local systems of $\G_m$. Hence,
we necessarily should have
$$(-1)^kK_{\psi,\chi_1,\ldots,\chi_k}=
c(-1)^lt_{a^{-1}}^*K_{\psi,\eta_1,\ldots,\eta_l}.$$
Making the multiplicative Fourier transform we conclude that the
equality
$$
\prod_{i=1}^k (-g(\la\chi_i))=c\cdot\la(a)\cdot
\prod_{j=1}^l (-g(\la\eta_j))
$$
holds {\it for all} $\la\in X(\F_q^*)$.
Now considering the jumps of the
absolute value of both sides
we immediately derive that the sets $\{\chi_1,\ldots,\chi_k\}$
and $\{\eta_1,\ldots,\eta_l\}$ are the same.
\ed

\section{Identities with the Fourier transform}

\subsection{Main theorem}

Let $\la_1,\ldots,\la_k$ be characters of
$\F_q^*$, $(n_1,\ldots,n_k)$ be a collection of integers such that 
$(n_i,q)=1$ for every $i$. Let us denote
\begin{equation}\label{Inla}
I^{n_1,\ldots,n_k}_{\la_1,\ldots,\la_k}(a)=
\sum_{(x_1,\ldots,x_k)\in(\F_q^*)^k}\psi(a\prod_{i=1}^kx_i^{n_i})
\la_1(x_1)\ldots\la_k(x_k)
\end{equation}
where $a\in\F_q^*$.

\begin{lem}\label{roots} Let $\la$ be a character of $\F_q^*$.
Then for any $d>0$ such that $(d,q)=1$ and any $a\in\F_q^*$, one has
$I^d_{\la}(a)=0$ unless there exists a character
$\mu$ such that $\la=\mu^d$. On the other hand,
$$I^d_{\mu^d}(a)=\sum_{\chi:\chi^d=1} g(\mu\chi)(\mu\chi)(a^{-1}).$$
\end{lem} 

\Pf. If $\la$ is not of the form $\mu^d$ then the restriction of
$\la$ to the subgroup of roots of unity of $d$-th order is non-trivial.
Thus, summing over cosets of this subgroup we get $I^d_{\la}(a)=0$.
Let $[d]:\F_q^*\ra\F_q^*$ be the homomorphism of raising to the $d$-th
power. Then we have
$$\sum_{\chi^d=1}\chi(x)=\cases 0 & x\not\in [d](\F_q^*),\\
                                d_1 & x\in [d](\F_q^*) \endcases   
$$
where $d_1=gcd(d,q-1)$.
Hence, 
$$\sum_{t\in\F_q^*}\psi(at^d)\mu(t^d)=d_1\cdot
\sum_{x\in [d](\F_q^*)}\psi(ax)\mu(x)=
\sum_{x\in\F_q^*,\chi^d=1}\psi(ax)(\mu\chi)(x).$$
\ed

\begin{lem}\label{int} 
One has $I^{n_1,\ldots,n_k}_{\la_1,\ldots,\la_k}(a)=0$
unless there exists $\la$ such that $\la_i=\la^{n_i}$ for all $i$.
One has 
$$I^{n_1,\ldots,n_k}_{\la^{n_1},\ldots,\la^{n_k}}(a)=
(q-1)^{k-1}\sum_{t\in\F_q^*}\psi(a t^d)\la^d(t)=(q-1)^{k-1}
\sum_{\chi:\chi^d=1} g(\la\chi)(\la\chi)(a^{-1}).$$
where $d=gcd(n_1,\ldots,n_k)$.
\end{lem}

\Pf . Let $n_i=n_i'd$. Since $gcd(n'_1,\ldots,n'_k)=1$ we can
choose new coordinates $y_i=\prod_j x_j^{a_{ij}}$ on $\G_m^k$
such that $y_1=\prod x_j^{n'_j}$.
Then we have
$$I^{n_1,\ldots,n_k}_{\la_1,\ldots,\la_k}(a)=
\sum_{(y_1,\ldots,y_k)\in(\F_q^*)^k} \psi(ay_1^d)\prod_{i,j}
\la_i(y_j^{b_{ij}})$$ 
where $(b_{ij})$ is the inverse matrix to $(a_{ij})$. Therefore,
we get zero unless $\prod_i\la_i^{b_{ij}}=1$ for every $j>1$,
i.e. 
$$\prod_{i,j}\la_i^{b_{ij}}(y_j)=\la'(y_1)$$
where $\la'=\prod_i \la_i^{b_{i1}}$. But the LHS is just
$\prod_i \la_i(x_i)$. Thus, the condition is that
$$\prod_i \la_i(x_i)=\la'(\prod_i x_i^{n'_i}),$$ 
i.e. $\la_i=(\la')^{n'_i}$. It remains to apply Lemma \ref{roots}.
\ed

Let us denote by 
$F^{n_1,\ldots,n_k}_{\chi_1,\ldots,\chi_k}(a)=
F^{n_1,\ldots,n_k}_{\chi_1,\ldots,\chi_k}(a,\psi)$ the
simple perverse sheaf on $\A^k$ obtained as the Goreski-MacPherson
extension of the smooth perverse sheaf 
$$L_{\psi}(a\prod_i x_i^{n_i})\otimes 
L_{\chi_1}\boxtimes\ldots\boxtimes L_{\chi_k}[k]$$
on $\G_m^k$. For $a=1$ we get the sheaf which we earlier denoted
$F^{n_1,\ldots,n_k}_{\chi_1,\ldots,\chi_k}$.

Recall that for every character $\chi\in X(\F_q^*)$ and $n\neq 0$ 
we have defined the
divisor $D_{\chi,n}\in\Div(X)$ by (\ref{Dchin}) and (\ref{Dchin2}).

\begin{thm}\label{monom}
Assume that $(n_i,q)=1$ for all $i$. One has an isomorphism
$$\FF(F^{n_1,\ldots,n_k}_{\chi_1,\ldots,\chi_k}(a))\simeq
V\otimes F^{m_1,\ldots,m_k}_{\eta_1,\ldots,\eta_k}(b)$$ 
where $V$ is a one-dimensional $\ov{\Q}_l$-vector space with
$\Gal(\ov{\F_q}/\F_q)$-action in the following situations:

\noindent
(i) $\sum n_i=2$, $m_i=n_i$ for all $i$;
$\eta_i=\chi^{n_i}\chi_i^{-1}$
where the characters $\chi,\chi_i\in X(\F_q^*)$ satisfy
\begin{equation}\label{div1}
D_{1,1}+D_{\chi^{-1},1}=\sum_i D_{\chi_i^{-1},n_i};
\end{equation}
if $gcd(n_1,\ldots,n_k)=2$ then we require that $\chi$
is the non-trivial character of order $2$ (so $q$ should be odd);
\begin{equation}\label{ab}
ab=-\prod_i n_i^{-n_i};
\end{equation}
\begin{equation}\label{V1}
V=G(\chi^{-1})\otimes(\bigotimes_i G(\chi_i))\otimes
(L_{\chi})_{-b}(-m)
\end{equation}
where $G(\la)$ are defined by (\ref{Gla}) and (\ref{G1}),
$F\mapsto F(1)$ is the Tate twist (the action of $\Frob_q$ is
multiplied by $q^{-1}$), $2m+1$ is the number of trivial characters among
$\chi,\chi_1,\ldots,\chi_k$.

\noindent
(ii) $\sum n_i=0$, $m_i=-n_i$ for all $i$;
$\eta_i=\chi^{n_i}\chi_i^{-1}$
where the characters $\chi,\chi_i\in X(\F_q^*)$ satisfy
\begin{equation}\label{div2}
D_{1,1}-D_{\chi^{-1},1}=\sum_i D_{\chi_i^{-1},n_i};
\end{equation}
if $gcd(n_1,\ldots,n_k)>1$ then we require $\chi=1$;
\begin{equation}\label{a/b}
\frac{a}{b}=\prod_i n_i^{-n_i};
\end{equation}
\begin{equation}\label{V2}
V=G(\chi)\otimes(\bigotimes_i G(\chi_i))\otimes
(L_{\chi^{-1}})_{-b}(-m)
\end{equation}
where $2m+1$ is the number of trivial characters among
$\chi,\chi_1,\ldots,\chi_k$.

\end{thm}

\Pf . Note that for any field extension $\F_q\sub\F_{q_1}$ the extension
of scalars of $F^{n_1,\ldots,n_k}_{\chi_1,\ldots,\chi_k}(a,\psi)$ to
$\F_{q_1}$ is isomorphic to 
$F^{n_1,\ldots,n_k}_{\chi_1\circ\Nm,\ldots,\chi_k\circ\Nm}(a,\psi\circ\Tr)$. 
Notice also that our assumptions
do not change after arbitrary extension of scalars. Thus, by Lemma 
\ref{mainlemma} it suffices to prove the identity 
$$(-q)^k(\psi(a\prod_ix_i^{n_i}\chi_i(x_i)),\la_1\otimes\ldots\otimes\la_k)=
c\cdot\prod_{i=1}^k \ov{g(\la_i)} 
\cdot(\psi(b\prod_ix_i^{m_i}\eta_i(x_i)),
\la_1^{-1}\otimes\ldots\otimes\la_k^{-1})$$
for every collection of non-trivial characters $\la_i$,
where $c=\Tr(\Frob_q,V)$.
Using the previous notation we can write this identity as follows:
\begin{equation}\label{auxid}
(-q)^k I^{n_1,\ldots,n_k}_{\frac{\chi_1}{\la_1},\ldots,
\frac{\chi_k}{\la_k}}(a)=c\cdot
\prod_{i=1}^k \ov{g(\la_i)} \cdot
I^{m_1,\ldots,m_k}_{\eta_1\la_1,\ldots,\eta_k\la_k}(b).
\end{equation}
Now let us specialize to different cases.

\noindent
(i) In this situation both sides are zero unless there exists $\la$
such that 
$$\frac{\chi_i}{\la_i}=\la^{n_i}.$$ 
Then we have $\la_i=\chi_i\la^{-n_i}$, so that
$\eta_i\la_i=(\chi\la^{-1})^{n_i}$.
Assume first
that $gcd(n_1,\ldots,n_k)=1$. Then according to Lemma \ref{int}
our identity (\ref{auxid}) takes form
$$(-q)^k g(\la)\la(a^{-1})=c\cdot\prod_{i=1}^k 
\ov{g(\chi_i\la^{-n_i})} \cdot g(\chi\la^{-1})
\cdot(\chi\la^{-1})(b^{-1}). 
$$
Notice that the relation (\ref{div1}) implies that there exists $i$
such that $\chi_i=1$ and $j$ such that $\chi_j=\chi^{n_j}$.
Therefore, the non-triviality of all the characters $\la_i=\chi_i\la^{-n_i}$
implies the non-triviality of $\la$ and of $\chi\la^{-1}$.
On the other hand, we have $|c|=q^{\frac{k}{2}}$. Therefore, both sides
of our identity have the same absolute value so it suffices to prove the
identity modulo $q^{\Z}$. Then we can rewrite it
as follows:
$$(-1)^k g(\la)g(\la\chi^{-1})
\equiv c\cdot(\prod_i \chi_i)(-1)\cdot
(\chi)(-b^{-1})\la(-ab)\cdot\prod_{i=1}^k
g(\chi_i^{-1}\la^{n_i})\mod q^{\Z}.
$$
On the other hand, by Corollary \ref{maincor} and (\ref{div1})
we have
$$g(\la)g(\la\chi)\equiv - g(\chi^{-1})
\prod_i\frac{g(\chi_i^{-1}\la^{n_i})}{\la(n_i^{n_i})g(\chi_i^{-1})}\mod
q^{\Z}.$$
Substituting this in the previous identity and using (\ref{ab})
we get
$$c\equiv (-1)^{k+1}\chi(-b)g(\chi^{-1})\prod_{i=1}^k g(\chi_i)
\mod q^{\Z}$$
which follows from (\ref{V1}).

Now consider the case $gcd(n_1,\ldots,n_k)=2$.
Then using Lemma \ref{int} and the equality $s=N/2$
the identity (\ref{auxid}) can be rewritten as
\begin{align*}
&(-q)^k \left(g(\la)\la(a^{-1})+g(\la\eps_2)(\la\eps_2)(a^{-1})
\right)=\\
&c\cdot
\prod_{i=1}^k \ov{g(\chi_i\la^{-n_i})}\cdot 
\left(g(\eps_2\la^{-1})(\eps_2\la^{-1})(b^{-1})+
g(\la^{-1})\la^{-1}(b^{-1})\right)
\end{align*}
where $\eps_2$ is the non-trivial character of order $2$.
Since the characters $\la$ and $\la\eps_2$ are non-trivial
we have $g(\la^{-1})=q\la(-1)g(\la)^{-1}$ and
$g(\la^{-1}\eps_2)=q\la(-1)\eps_2(-1)g(\la\eps_2)^{-1}$. Notice also
that since $n_i$ are even we have $\eps_2(-ab)=1$. Therefore, our
identity follows from 
$$(-q)^kg(\la)g(\la\eps_2)=q\la(-ab)\eps_2(-b^{-1})c\cdot
\prod_{i=1}^k \ov{g(\chi_i\la^{-n_i})},$$
which can be proven as in the previous case.

\noindent
(ii) Again both sides of (\ref{auxid}) are zero unless there exists $\la$
such that $\la_i=\chi_i\la^{-n_i}$ so that
$\eta_i\la_i=(\chi\la^{-1})^{n_i}=(\chi^{-1}\la)^{m_i}$.
Assume first that $gcd(n_1,\ldots,n_k)=1$.
Then according to Lemma \ref{int}
the identity (\ref{auxid}) takes form
\begin{equation}\label{auxid2}
(-q)^k g(\la)\la(a^{-1})=c\cdot\prod_{i=1}^k 
\ov{g(\chi_i\la^{-n_i})} \cdot g(\chi^{-1}\la)
\cdot(\chi^{-1}\la)(b^{-1}). 
\end{equation}
In the case $gcd(n_1,\ldots,n_k)=d>1$
the identity (\ref{auxid}) is equivalent to
\begin{align*}
(-q)^k \sum_{\eta^d=1} g(\la\eta)(\la\eta)(a^{-1})=c\cdot\prod_{i=1}^k 
\ov{g(\chi_i\la^{-n_i})} \cdot 
\sum_{\eta^d=1} g(\la\eta)(\la\eta)(b^{-1}). 
\end{align*}
which reduces to (\ref{auxid2}) with $\chi=1$ 
since $\eta(a)=\eta(b)$ for any character
$\eta$ of order $d$ (note that by (\ref{a/b}) the ratio $a/b$ is the 
$d$-th power). 
Note that $|c|=q^{\frac{k}{2}}$. On the other hand, if $\chi\neq 1$
then the non-triviality of characters $\chi_i^{-1}\la^{n_i}$ and
the equality (\ref{div2}) imply the non-triviality of $\chi$ and
$\chi^{-1}\la$. Thus, both sides of (\ref{auxid2})
have the same absolute value so we can work modulo $q^{\Z}$.
The remaining part of the proof is similar to the case (i).
\ed

\subsection{Hypergeometric sheaves}

If we work over the algebraic closure of a finite field then
the isomorphisms of Theorem (\ref{monom}) follow easily from
the theory of hypergeometric sheaves developed in \cite{K1}, \cite{K2}
and \cite{GL}.
Indeed, let $f:\G_m^n\ra\G_m$ be a non-constant homomorphism of tori.
Then for any character $\chi$ of $\G_m^n(\F_q)$ we can consider
the sheaf $j_{!*}L(\psi(f(x))\chi(x))[n]$ on $\A^n$ where
$j:\G_m^n\ra\A^n$ is the standard open embedding. By a simple coordinate
change one can see that 
$$j^*\FF j_{!*}L(\psi(f(x))\chi(x))[n]\simeq (f^{-1})^*j_1^*\FF j_{1,!*}H$$
where $j_1:\G_m\ra\A^1$ is the standard embedding,
$H$ is the hypergeometric sheaf on $\G_m$ defined as follows:
$$H=\im(f_!L(\psi(\sum_i x_i)\chi(x))[n]\ra
f_*L(\psi(\sum_i x_i)\chi(x))[n]).$$
To see that $H$ is a hypergeometric and to compute it explicitly
we notice that according to Proposition 5.6.2 of \cite{K1}
for any multiplicative character $\eta$ one has
$$[N]_*L(\psi(t)\eta(t))[1]\simeq\Hyp(!,\psi;(\eta_i);\emptyset)$$
where $[N]:\G_m\ra\G_m$ is the morphism of
raising to the $N$-th power ($N$ is assumed to be relatively prime to $q$),
$(\eta_i)$ is set of $N$-th roots of $\eta$ considered as a character of
$\ov{\F_q}^*$, $\Hyp(!,\psi;(\eta_i);\emptyset)$ is the hypergeometric
sheaf defined as the $!$-convolution of the $N$ sheaves
$L(\psi(t)\eta_i(t))[1]$
on $\G_m$. On the other hand, the Cancellation Theorem 8.4.7 of \cite{K2}
implies that for arbitrary collections of characters
$(\eta_i)$ and $(\rho_j)$ the unique simple quotient of the
$!$-convolution of the sheaves $L(\psi(t)\eta_i(t))[1]$ and
$L(\psi(-t^{-1})\rho_j)[1]$ on $\G_m$ depends only on the
divisor $\sum_i (\eta_i) -\sum_j (\rho_j)$. According to Proposition 8.1.4
of \cite{GL} this unique simple quotient coincides with the image of the
natural morphism from the $!$-convolution to the $*$-convolution of the
same collection of sheaves. Thus, $H$ is the (irreducible)
hypergeometric sheaf corresponding to the divisor $\sum D_{\chi_i,n_i}$
where $f(x)=\prod_i x_i^{n_i}$, $\chi(x)=\prod_i\chi_i(x_i)$.
It remains to notice that
$j_1^*\FF(j_{!*}H)$ is the image of the natural morphism
$L_{\psi}*_!\inv^*H\ra L_{\psi}*_*\inv^*H$, where $\inv:\G_m\ra\G_m$
is the inversion. Thus, it is also a hypergeometric
sheaf which can be computed by Cancellation Theorem. This gives
isomorphisms of Theorem \ref{monom} over an algebraically closed field.
On the other hand, this argument allows to compute the rank of the Fourier
transform of the sheaf
$$j_{!*}L(\psi(\frac{\prod_i x_i^{m_i}}{\prod_j y_j^{n_j}})\prod_i\chi_i(x_i)
\prod_j \eta_j(y_j))$$
for generic characters $\chi_i$, $\eta_j$
(i.e. when there is no cancellation), where $m_i$ and $n_j$ are positive
and prime to $q$. Namely, using Theorem 8.4.2 of \cite{K2} we find that
this rank is equal to max$(\sum_i m_i, \sum_j n_j+1)$.

\subsection{Examples}\label{examples}

In the case $k=1$ the conditions
of Theorem \ref{monom} are satisfied only in the case $n_1=N=2$,
$\chi_1=1$, $\chi=\eps_2$, which corresponds to the isomorphism
$$\FF(L_{\psi}(ax^2))\simeq G(\eps_2)\otimes (L_{\eps_2})_{a})
\otimes L_{\psi}(-\frac{x^2}{4a}).$$

For $k=2$ we have the following examples:

\noindent
1. $n_1=n_2=1$, $\chi_1=\eta_2=1$, $\chi_2=\eta_1=\chi$,
where $\chi$ is an arbitrary character, the corresponding
isomorphism is
$$\FF(L(\psi(axy)\chi(x))\simeq (L_{\chi})_{-a^{-1}}(-1)\otimes
L(\psi(-a^{-1}xy)\chi(y)).$$

\noindent
2. $n_1=3$, $n_2=-1$, $\chi_1=\eta_1=1$, $\chi_2=\eta_2=\eps_3$
where $\eps_3$ is a non-trivial character of order $3$ (so $q-1$
is divisible by $3$). In this
example we have
$$\FF(j_{!*}L(\psi(\frac{x^3}{y})\eps_3(y))[2])\simeq \ov{\Q}_l(-1)\otimes
j_{!*}L(\psi(\frac{x^3}{27y})\eps_3(y))[2].$$
Let $f_{3,-1}:\F_q^2\ra\C$ be the trace function of
$F^{3,-1}_{1,\eps_3}$. Then according to Theorem \ref{2dim} we have
$$f_{3,-1}(x,y)=\cases \psi(\frac{x^3}{y})\eps_3(y), & y\neq 0,\\
                0, & y=0, x\neq 0,\\
                g(\eps_3^{-1}), & x=y=0.\endcases
$$
Now Theorem \ref{monom} implies that
$$\widehat{f_{3,-1}}(x,y)=qf_{3,-1}(x,27y).$$

\noindent
3. $n_1=4$, $n_2=-2$, $\chi_1=\eta_1=1$, $\chi_2=\eta_2=\eps_2$ (so
$q$ is odd). We have
$$\FF(j_{!*}L(\psi(a\frac{x^4}{y^2})\eps_2(y))[2])\simeq\ov{\Q}_l(-1)
\otimes (L_{\eps_2})_a
\otimes j_{!*}L(\psi(\frac{x^4}{2^{10}ay^2})\eps_2(y))[2].$$
Let $f^a_{4,-2}:\F_q^2\ra\C$ be the trace function of
$F^{4,-2}_{1,\eps_2}(a)$. Then using Theorem \ref{2dim} and Lemma
\ref{roots} we find
$$f^a_{4,-2}(x,y)=\cases \psi(a\frac{x^4}{y^2})\eps_2(y), & y\neq 0,\\
                0, & y=0, x\neq 0,\\
                g(\eps_4)\eps_4(a^{-1})+g(\eps_4^{-1})\eps_4(a), & x=y=0
                \endcases
$$
if $q\equiv 1\mod(4)$. In case $q\equiv 3\mod(4)$ the function
$f^a_{4,-2}(x,y)$ is just an extension by zero of
$\psi(a\frac{x^4}{y^2})\eps_2(y)$. Now
Theorem \ref{monom} implies that
$$\widehat{f^a_{4,-2}}(x,y)=q\eps_2(a)f^a_{4,-2}(x,32y).$$

\noindent
4. $n_1=-n_2=n>0$, $\chi_1=\chi_2^{-1}=\eta_1^{-1}=\eta_2=\chi$,
where $\chi$ is any character. We have
$$\FF(j_{!*}L(\psi(a\frac{x^n}{y^n})\chi(\frac{x}{y}))[2])
\simeq\ov{\Q}_l(-1)
\otimes j_{!*}L(\psi((-1)^na\frac{y^n}{x^n})\chi(-\frac{y}{x}))[2].$$
If $n=1$ there is one more example:
$\chi_1=\eta_2=1$, $\chi_2=\eta_1^{-1}=\chi$ is any character,
the corresponding isomorphism is
$$\FF(j_{!*}L(\psi(\frac{x}{y})\chi(y))[2])
\simeq\ov{\Q}_l(-1)
\otimes j_{!*}L(\psi(-\frac{y}{x})\chi^{-1}(-x))[2].$$
Again using Theorem \ref{2dim} we can extract the identities
for Fourier transforms of functions of $\F_q^2$ from these isomorphisms.
In particular, we obtaing the following identity (which we'll use later):
\begin{equation}\label{psixy}
\sum_{(x,y)\in(\F_q^*)^2}\psi(a\frac{x}{y}+x\hat{x}+y\hat{y})
\chi(\frac{x}{y})=q\psi(-a\frac{\hat{y}}{\hat{x}})\chi(-\frac{\hat{y}}
{\hat{x}})-g(\chi)\chi^{-1}(a)
\end{equation}
for any non-zero $a$, $\hat{x}$ and $\hat{y}$ and any character $\chi$
(note that this identity is also easy to check directly by changing
variables to $t=x/y$ and $y$).

Now let us consider some higher-dimensional examples.

\noindent
5. If the conditions of Theorem \ref{monom} are satisfied for
the collection $(n_1,\ldots,n_k)$, $(\chi_1,\ldots,\chi_k)$,
then they are also satisfied for the collection
$(n_1,\ldots,n_k,1,-1)$, $(\chi_1,\ldots,\chi_k,\chi,\chi^{-1})$
where $\chi$ is any character.
For example,
for any characters $\chi_1,\ldots,\chi_{k+1}$ we have
\begin{align*}
&\FF(j_{!*}L(\psi(\frac{x_1\ldots x_{k+2}}{y_1\ldots y_k})
\chi_1(\frac{x_1}{y_1})\ldots\chi_k(\frac{x_k}{y_k})\chi_{k+1}(x_{k+1}))
[2k+2])\simeq\ov{\Q}_l(-k-1)\otimes\\ 
&j_{!*}L(\psi((-1)^{k+1}\frac{x_1\ldots x_{k+2}}{y_1\ldots y_k})
\frac{\chi_{k+1}}{\chi_1}(-\frac{x_1}{y_1})\ldots
\frac{\chi_{k+1}}{\chi_k}(-\frac{x_k}{y_k})\chi_{k+1}(-x_{k+2}))[2k+2],
\end{align*}

\begin{align*}
&\FF(j_{!*}L(\psi(\frac{x_1\ldots x_{k}}{y_1\ldots y_k})
\chi_1(\frac{x_1}{y_1})\ldots\chi_k(\frac{x_k}{y_k}))
[2k])\simeq\\
&\ov{\Q}_l(-k)\otimes
j_{!*}L(\psi((-1)^k\frac{y_1\ldots y_{k}}{x_1\ldots x_k})
\chi_1(-\frac{y_1}{x_1})\ldots
\chi_k(-\frac{y_k}{x_k})[2k],
\end{align*}

\begin{align*}
&\FF(j_{!*}L(\psi(\frac{x_1\ldots x_{k+1}}{y_1\ldots y_{k+1}})
\chi_1(\frac{x_1}{y_1})\ldots\chi_k(\frac{x_k}{y_k}))\chi_{k+1}(y_{k+1})
[2k+2])\simeq\ov{\Q}_l(-k-1)\otimes\\
&(L_{\prod_{i=1}^k\chi_i})_{-1}\otimes
j_{!*}L(\psi((-1)^{k+1}\frac{y_1\ldots y_{k+1}}{x_1\ldots x_{k+1}})
\chi_1(\frac{y_1}{x_1})\ldots
\chi_k(\frac{y_k}{x_k})\chi_{k+1}^{-1}((-1)^{k+1}x_{k+1}))[2k+2],
\end{align*}

\noindent
6. More generally, if the conditions of Theorem \ref{monom} are satisfied for
the expression $\psi(\prod_i x_i^{n_i})\prod_i\chi_i(x_i)$,
then they are also satisfied for the expressions 
$$\psi(\prod_i x_i^{n_i}\cdot \frac{u^k}{v_1\ldots v_k})
\prod_i \chi_i(x_i)\cdot\frac{\prod_{j=1}^k(\chi\eps_k^j)(v_j)}{\chi^k(u)} 
$$
and
$$\psi(\prod_i x_i^{n_i}\cdot \frac{v_1\ldots v_k}{u^k})
\prod_i \chi_i(x_i)\cdot
\frac{\prod_{j=1}^k(\chi\eps_k^j)(v_j)}{\chi^k(u)} 
$$
where $\chi$ is any character.

For example, for any $k>0$ and any characters $\chi,\eta$ we have
\begin{align*}
&\FF(j_{!*}L(\psi(\frac{x_1\ldots x_{k+2}}{y^k})
\frac{\prod_{i=1}^k
(\chi\eps_k^i)(x_i)}{\chi^k(y)}\cdot\eta(x_{k+1}))[k+3])\simeq\ov{\Q}_l(-2)
\otimes\\
&(\bigotimes_{i=1}^{k-1}G(\eps_k^i))
j_{!*}L(\psi((-1)^{k+1}k^k\cdot\frac{x_1\ldots x_{k+2}}{y^k})
\frac{\prod_{i=1}^k(\chi^{-1}\eta\eps_k^{-i})(x_i)}
{(\chi^{-1}\eta)^k(-y/k)}\cdot\eta(-x_{k+2}))[k+3],
\end{align*}

\begin{align*}
&\FF(j_{!*}L(\psi(\frac{x_1\ldots x_k}{y^k})
\frac{\prod_{i=1}^k(\chi\eps_k^i)(x_i)}{\chi^k(y)})[k+1])\simeq\\
&\ov{\Q}_l(-1)\otimes(\bigotimes_{i=1}^{k-1}G(\eps_k^i))
\otimes 
j_{!*}L(\psi((-1)^k\frac{y^k}{k^kx_1\ldots x_k})
\frac{\prod_{i=1}^k(\chi^{-1}\eps_k^{-i})(x_i)}{\chi^{-k}(-y/k)})[k+1].
\end{align*}

\noindent
7. Example 2 has the following generalization: 
for any $k>0$ and any $1\le i\le k+1$ one has
\begin{align*}
&\FF(j_{!*}L(\psi(\frac{x^{k+2}}{y_1\ldots y_k})\prod_{1\le j<i}
\eps_{k+2}^j(y_j)\prod_{i\le j\le k}\eps_{k+2}^{j+1}(y_j))[k+1])\simeq\\
&\bigotimes_{j=1}^{k+1}G(\eps_{k+2}^j)
\otimes j_{!*}L(\psi((-1)^{k+1}\frac{x^{k+2}}{(k+2)^{k+2}y_1\ldots y_k})
\prod_{1\le j<i}
\eps_{k+2}^{i-j}(y_j)\prod_{i\le j\le k}\eps_{k+2}^{i-j-1}(y_j))[k+1]
\end{align*}

\noindent
8. Here is an example involving the monom $\frac{x_1\ldots x_k}{y^k}$ which
is different from the one considered in example 6. Let $\chi$
be any character, then for any $i$, $1\le i\le k$, we have 
\begin{align*}
&\FF(j_{!*}L(\psi(\frac{x_1\ldots x_k}{y^k})
\frac{\prod_{1\le j<i}(\chi\eps_k^j)(x_j)\prod_{i\le j\le k-1}
(\chi\eps_k^{j+1})(x_j)}{\chi^k(y)})[k+1])\simeq\\
&\ov{\Q}_l(-1)\otimes(\bigotimes_{j=1}^{k-1}G(\eps_k^j))\otimes
\\
& j_{!*}L(\psi((-1)^k\frac{y^k}{k^kx_1\ldots x_k})
\frac{\prod_{1\le j<i}\eps_k^{i-j}(x_j)\prod_{i\le j\le k-1}
\eps_k^{i-j-1}(x_j)}{(\chi^{k-1}\eps_k^{-i})(y)}(\chi\eps_k^i)(-x_k))[k+1])
\end{align*}

\subsection{Identities with binomial coefficients}

Combining the example 5 of the previous section with Theorem \ref{GMtr}
we obtain some identities with polynomials (\ref{anm}) and (\ref{bnm}).
We need the  following simple lemma.

\begin{lem}\label{rs}
For any character $\chi$ of $\F_q^*$ and any
$a,\hat{x}_1,\ldots,\hat{x}_n,\hat{y}_1,\ldots,\hat{y}_n\in\F_q^*$
one has
\begin{align*}
&\sum_{(x,y)\in(\F_q^*)^{2n}}
\psi(\frac{x_1\ldots x_n}{y_1\ldots y_n}+
\sum_{m=1}^n(x_m\hat{x}_m+y_m\hat{y}_m))
\chi(\frac{x_1\ldots x_n}{y_1\ldots y_n})=\\
&q^n
\psi((-1)^n\frac{\hat{y}_1\ldots \hat{y}_n}{\hat{x}_1\ldots \hat{x}_n})
\chi((-1)^n\frac{\hat{y}_1\ldots \hat{y}_n}{\hat{x}_1\ldots \hat{x}_n})
-\frac{q^n-1}{q-1}g(\chi).
\end{align*}
\end{lem}

\Pf . Use the equality (\ref{psixy}) and induction in $n$.
\ed

Consider the isomorphism
$$
\FF(j_{!*}L(\psi(\frac{x_1\ldots x_n}{y_1\ldots y_n})))\simeq
\ov{\Q}_l(-n)\otimes
j_{!*}
L(\psi((-1)^n\frac{\hat{y}_1\ldots \hat{y}_n}{\hat{x}_1\ldots \hat{x}_n})).
$$
Let us restrict it to the point with
$\hat{y}_1=\ldots=\hat{y}_r=\hat{x}_1=\ldots=\hat{x}_s=0$ and
$\hat{y}_{r+1}\ldots \hat{y}_n\hat{x}_{s+1}\ldots\hat{x}_n\neq 0$.
Taking the traces of $\Frob_q$ we get the identity
\begin{align*}
&\sum_{(x,y)\in(\F_q^*)^{2n}}
\psi(\frac{x_1\ldots x_n}{y_1\ldots y_n}+
\sum_{m=1}^n(x_m\hat{x}_m+y_m\hat{y}_m))+\\
&\sum_{(i,j,k,l)\neq(0,0,0,0)}N(i,j,k,l)a(i+j,k+l)=\cases
q^na(r,s), & (r,s)\neq(0,0)\\
q^n\psi((-1)^n\frac{\hat{y}_1\ldots \hat{y}_n}{\hat{x}_1\ldots \hat{x}_n}),
& r=s=0 \endcases
\end{align*}
where
$$N(i,j,k,l)=\sum_{(x,y)\in S(i,j,k,l)}
\psi(\sum_{m=1}^n (x_m\hat{x}_m+y_m\hat{y}_m))$$
and $S(i,j,k,i)$ is the set of $(x,y)\in\F_q^n\times\F_q^n$
with exactly $i$ coordinates $(x_1,\ldots x_r)$ vanish,
$j$ coordinates $(x_{r+1},\ldots,x_n)$ vanish,
$k$ coordinates $(y_{1},\ldots,y_s)$ vanish, and
$l$ coordinates $(y_{s+1},\ldots,y_n)$ vanish.
Thus,
$$N(i,j,k,l)={r\choose i}{n-r\choose j}{s\choose k}{n-s\choose l}
(q-1)^{r+s-i-k}(-1)^{r+s+j+l}.$$
On the other hand, it is easy to see that for $(r,s)\neq (0,0)$ we have
$$\sum_{(x,y)\in(\F_q^*)^{2n}}
\psi(\frac{x_1\ldots x_n}{y_1\ldots y_n}+
\sum_{m=1}^n(x_m\hat{x}_m+y_m\hat{y}_m))=(q-1)^{r+s-1}(-1)^{r+s+1}$$
Thus, we arrive to the following identity for any $(r,s)\neq(0,0)$:
\begin{equation}\label{binom1}
\begin{array}{l}
\sum_{(i,j,k,l)\neq(0,0,0,0)}
{r\choose i}{n-r\choose j}{s\choose k}{n-s\choose l}
(q-1)^{r+s-i-k}(-1)^{r+s+j+l}a(i+j,k+l)=\\
q^na(r,s)+(q-1)^{r+s-1}(-1)^{r+s}.
\end{array}
\end{equation}
In the case $r=s=0$ using Lemma \ref{rs} we get the identity
\begin{equation}\label{binom2}
\sum_{(j,l)\neq(0,0)}
{n\choose j}{n\choose l}(-1)^{j+l}a(j,l)=-\frac{q^n-1}{q-1}.
\end{equation}

Similarly the isomorphism
$$\FF(j_{!*}L(\psi(\frac{x_1\ldots x_n}{y_1\ldots y_n})
\chi(\frac{x_1\ldots x_n}{y_1\ldots y_n})))\simeq\ov{\Q}_l(-n)\otimes
j_{!*}L(\psi((-1)^n\frac{y_1\ldots y_n}{x_1\ldots x_n})
\chi((-1)^n\frac{y_1\ldots y_n}{x_1\ldots x_n}))
$$
for a non-trivial character $\chi$ leads to the identities
\begin{equation}\label{binom3}
\begin{array}{l}
\sum_{(i,j,k,l)\neq(0,0,0,0)}
{r\choose i}{n-r\choose j}{s\choose k}{n-s\choose l}
(q-1)^{r+s-i-k}(-1)^{r+s+j+l}b(i+j,k+l)=\\
q^nb(r,s)+(q-1)^{r+s-1}(-1)^{r+s+1}
\end{array}
\end{equation}
for $(r,s)\neq (0,0)$ and
\begin{equation}\label{binom4}
\sum_{(j,l)\neq(0,0)}
{n\choose j}{n\choose l}(-1)^{j+l}b(j,l)=\frac{q^n-1}{q-1}.
\end{equation}

%These identities are not difficult to check directly 

%\begin{align*}
%&\FF(j_{!*}L_{\psi}(\frac{x_1\ldots x_{n+2}}{y_1\ldots y_n}))\simeq\\
%&\ov{\Q}_l(-n-1)\otimes
%j_{!*}L_{\psi}((-1)^{n+1}\frac{x_1\ldots x_{n+2}}{y_1\ldots y_n}).
%\end{align*}

\section{Identities containing norms}
\label{norm-sec}

\subsection{More identities with Gauss sums}

We want to generalize Corollary \ref{maincor} to include the
identities for Gauss sums containing norms. 
In other words, we want to employ systematically both Hasse-Davenport
identities (\ref{HasDavsim}) and (\ref{HasDav}).
Let us fix a finite field $\F_q$. For every $d>0$, a character
$\chi\in X(\F_{q^d}^*)$, and an integer $n>0$, $gcd(n,q)=1$, consider
the function on $X(\F_q^*)$ defined by
$$f_{\chi,n}(\la)=\frac{g((\la^n\circ\Nm_d)\chi)}{\la(n^{nd})g(\chi)}$$
For $d=1$ this definition coincides with the one we considered before.
As before the map $[\chi,n]\mapsto f_{\chi,n}$ extends
to a homomorphism from $A^{(p)}(X(\F_{q^d}^*))$ to $\CC(X(\F_q^*),\C^*)$
due to the identity (\ref{HasDav}).

Let us define the abelian group $A_q$ as the quotient
of $\oplus_d A^{(p)}(X(\F_{q^d}^*))$ by the subgroup generated
by the elements of the form $k[\chi,n]-[\chi\circ\Nm_k,n]$
for all $\chi\in X(\F_{q^d}^*)$, $n,k\in\Z_{>0}$.
The identity (\ref{HasDavsim}) implies that for
every $k>0$ we have
$$f_{\chi\circ\Nm_k,n}=f_{\chi,n}^k.$$
Hence, the map $[\chi,n]\mapsto f_{\chi,n}$ extends to a homomorphism
from $A_q$ to $\CC(X(\F_q^*),\C^*)$.

On the other hand, we can define a homomorphism
$$\b:A_q\ra\Div(X)$$
by sending $[\chi,n]$ to the divisor $d\cdot D_{\chi,n}$,
where $\chi\in X(\F_{q^d}^*)$, $D_{\chi,n}$ is defined
by (\ref{Dchin}).

\begin{thm} $\ker\b=0$.
\end{thm}

\Pf . By definition the homomorphism $\b$ can be factorized as follows:
$$\b:A_q\stackrel{\b'}{\ra}\lim_d A^{(p)}(X(\F_{q^d}^*))\wt{\ra}\Div(X)$$
where the last arrow is the isomorphism (\ref{isdiv}), $\b'$
is defined by the formula $\b'([\chi,n])=d[\chi,n]$ for
$\chi\in X(\F_{q^d}^*)$. Clearly, $\b'$ is an isomorphism
modulo torsion so it suffices to prove that the group
$A_q$ is torsion-free. Fix a prime number $l$. Assume that
we have $lx=0$ for some $x\in A_q$. We can write $x$ in
the form 
$$x=\sum_i \pm [\chi_i,n_i]$$
with some $\chi_i\in X(\F_{q^{d_i}}^*)$. Let us write the degrees $d_i$
in the form $d_i=d'_il^{s_i}$ where $gcd(d'_i,l)=1$.
Set $d=\prod_i d'_i$. Using the relations in $A_q$ we obtain that
$x'=dx$ has form
$$x'=\sum_j \pm[\chi_j,n_j]$$
with $\chi_j\in X(\F_{q^{dl^{s_j}}}^*)$.
Since $d$ is relatively prime to $l$ it suffices to prove that $x'=0$.

Let $A_q(d,l)$ be the quotient of $\oplus_i A^{(p)}(\F_{q^{dl^i}}^*)$
by the subgroup generated by the elements $l[\chi,n]-[\chi\circ\Nm_l,n]$.
Then we have a sequence of homomorphisms
$$A_q(d,l)\stackrel{\b'_{d,l}}{\ra}\lim_i A^{(p)}(X(\F_{q^{dl^i}}^*))\hookrightarrow
\Div(X)$$
where $\b'_{d,l}([\chi,n])=dl^i[\chi,n]$ for $\chi\in X(\F_{q^{dl^i}}^*)$,
the last arrow is an embedding induced by (\ref{isdiv}). Notice, that
by Lemma \ref{dirsum} the group $A_q(d,l)$ has no torsion. Therefore,
the homomorphism $\b'_{d,l}$ being an isomorphism modulo torsion should
be injective. Hence, the composed map 
$$\b_{d,l}:A_q(d,l)\ra\Div(X)$$
is injective. Since $\b_{d,l}$ is a composition of $\b$ and of the natural
homomorphism $A_q(d,l)\ra A_q$ we deduce that $A_q(d,l)$ is a subgroup in
$A_q$. Now we have $x'\in A_q(d,l)$ and $lx'=0$. Since $A_q(d,l)$ has
no torsion this implies that $x'=0$.
\ed

This theorem allows to write the identities between Gauss sums
containing norms in the following form.

\begin{cor}\label{maincornorm} 
Let $\chi_1,\ldots,\chi_k$ be the collection of characters, 
$\chi_i\in X(\F_{q^{d_i}}^*)$; $n_1,\ldots,n_k$ be the collection of
integers such that $gcd(n_i,q)=1$.
Assume that $\sum_{i=1}^k d_iD_{\chi_i,n_i}=0$
in $\Div(X)$.
Then for every character $\la$ of $\F_q^*$ one has
$$\prod_i\frac{g((\la^{n_i}\circ\Nm_{d_i})\chi_i)}{\la(n_i^{n_id_i})g(\chi_i)}=
q^{m(\la)}$$
for some $m(\la)\in\Z$.
In particular, if $(\la^{n_i}\circ\Nm_{d_i})\chi_i\neq 1$ for all $i$ then
$2m(\la)$ is the number of $i$ such that $\chi_i=1$.
\end{cor}

One can rewrite this corollary in yet another form. Namely, let
$k$ be a finite \'etale $\F_q$-algebra, so that
$k=\prod_i \F_{q^{d_i}}$.
Let $X(k^*)$ be the group of characters of $k^*$. 
For every $\chi\in X(k^*)$ let us denote 
\begin{equation}\label{gk}
g(\chi)=\sum_{x\in k^*}\chi(x)\psi(\Tr_{k/\F_q}(x)).
\end{equation}
In fact, if $\chi=\prod_i\chi_i$, where $\chi_i\in X(\F_{q^{d_i}}^*)$ then
$$g(\chi)=\prod_i g(\chi_i,\psi\circ\Tr_{d_i}).$$
Let $V$ be a virtual finite module over $k$,
then $V$ is determined by the collection of
integers $(n_i=\rk_i V)$ such that $V=\sum_i n_i[\F_{q^{d_i}}]$.
For every character $\chi=\prod_i\chi_i$ of $k^*$
let us denote
$$D_{\chi,V}=\sum_i d_i D_{\chi_i,\rk_i V}.$$
Also set
$$p(V)=\prod_i n_i^{n_id_i}.$$
On the other hand, we can associate to $V$ the homomorphism
$$\sideset{}{_V}\det=\sideset{}{_{V/\F_q}}\det:
k^*\ra\F_q^*:(x_i)\mapsto \prod_i\Nm_{d_i}(x_i)^{n_i}.$$
Then we have
$$g((\la\circ\sideset{}{_V}\det)\chi)=
\prod_i g((\la^{n_i}\circ\Nm_{d_i})\chi_i,\psi\circ\Tr_{d_i}).$$
Thus, Corollary \ref{maincornorm} leads to the following statement.

\begin{thm}\label{gaussnorm} 
Let $V$ be a virtual finite $k$-module such that
$gcd(p(V),q)=1$, $\chi$ be
a character of $k^*$. Assume that $D_{\chi,V}=0$. Then for every
$\la\in X(\F_q^*)$ one has 
$$g((\la\circ\sideset{}{_V}\det)\chi)=q^{m(\la)}\cdot p(V)\cdot g(\chi)$$
where $m(\la)$ is an integer depending on $\la$.
\end{thm}

\begin{rem} In fact, our method allows to prove a stronger result.
Namely, consider the natural homomorphism $\Div(X)\ra\Div(X/\Frob_q)$.
Then one can replace the assumption $D_{\chi,V}=0$ by the weaker
assumption that the image of $D_{\chi,V}$ in $\Div(X/\Frob_q)$ is zero.
In this way one gets much more indentities between Gauss sums
(cf. \cite{De2},5.13). However, the importance of the identities of theorem   
\ref{gaussnorm} is that they hold universally over any finite extension of
a given finite field $\F_q$.
\end{rem}

\subsection{Identities with the Fourier transform containing norms}

Recall (see e.g. \cite{De1}) that for every finite \'etale $\F_q$-algebra $k$ 
one can define the ring scheme $\A^1k$ over $\F_q$ 
in a natural way, so that for every $\F_q$-algebra $A$ on has
$\A^1k(A)=k\otimes_{\F_q}A$.
Note that as an $\F_q$-scheme $\A^1k$ is non-canonically isomorphic to
$\A^d$, where $d=[k:\F_q]$. There is a natural
morphism of $\F_q$-schemes
$$\Tr:\A^1k\ra\A^1$$
inducing the usual trace map on points.
We also have an open subscheme $\G_mk\subset\A^1k$ of invertible elements
such that $\G_m k(\F_q)\simeq k^*$. For every character
$\chi\in X(k^*)$ we can define a rank 1 smooth l-adic sheaf
$L_{\chi}$ on $\G_m k$. Using the multiplication morphism
$$m:\A^1 k\times\A^1k\ra\A^1k$$
and the sheaf $\Tr^*L_{\psi}$ we can define the Fourier transform for
sheaves on $\A^1k$. 

Let us call a character
$\chi\in X(k^*)$ {\it non-degenerate} if $k=\prod_i\F_{q^{d_i}}$
and $\chi=\prod \chi_i$ where $\chi_i$ are non-trivial characters of
$\F_{q^{d_i}}$. 
For a non-degenerate character $\chi\in X(k^*)$ we have
$$\FF(j_!L_{\chi})\simeq G(\chi)\otimes j_!L_{\chi^{-1}}$$
where
\begin{equation}\label{Gk}
G(\chi)=H^d_c(\G_m k,L_{\chi}),
\end{equation}
is the one-dimensional $\ov{\Q}_l$-space on which
$\Frob_q$ acts as $(-1)^dg(\chi)$, where $d=[k:\F_q]$.
As before we will use the definition (\ref{Gk}) 
also in the case of degenerate characters.

For any extension $\F_q\sub\F_{q_1}$ we have
$$\A k\otimes_{\F_q}\F_{q_1}\simeq \A k'$$
where $k'=k\otimes_{\F_q}\F_{q_1}$, and for any $\chi\in X(k^*)$
$$L_{\chi}\otimes_{\F_q}\F_{q_1}\simeq L_{\chi\circ\Nm_{k'/k}}$$

Let $V$ be a virtual finite $k$-module such that
$gcd(p(V),q)=1$. We can define the scheme-theoretic version
of the homomorphism $\det_V$ considered above:
$$\sideset{}{_V}\det:\G_mk\ra\G_m.$$
Let $\chi$ be a character of $k^*$, $a$ be an element of $\F_q^*$.
We denote by
$F_{V,\chi}(a)=F_{V,\chi}(a,\psi)$ the
simple perverse sheaf on $\A^1k$
obtained as the Goreski-MacPherson
extension of the smooth perverse sheaf 
\begin{equation}\label{smoothperv}
L_{\psi}(a\sideset{}{_V}\det(x))\otimes L_{\chi}[d]
\end{equation}
on $\G_mk$, where $d=[k:\F_q]$.

We need a slight generalization of Lemma \ref{int}. Let us denote
$$I_{V,\la}(a)=\sum_{x\in k^*}\psi(a\sideset{}{_V}\det(x))\la(x)$$
where $V$ is a virtual finite $k$-module, $\la\in X(k^*)$, $a\in\F_q^*$.
For a virtual $k$-module $V$ with $\rk_i V=n_i$ let us denote 
$$d(V)=[F_q^*:\sideset{}{_V}\det(k^*)]=gcd(n_1,\ldots,n_r).$$

\begin{lem}\label{intnorm} Assume that $gcd(p(V),q)=1$.
One has $I_{V,\la}(a)=0$
unless there exists $\mu$ such that $\la=\mu\circ\det_V$.
One has 
$$I_{V,\mu\circ\sideset{}{_V}\det}(a)=
\frac{|k^*|}{q-1}\cdot
\sum_{\nu\in X(\F_q^*):\nu\circ\sideset{}{_V}\det=1} 
g(\mu\nu)(\mu\nu)(a^{-1}).$$
\end{lem}

For a virtual $k$-module $V=\sum_i n_i[\F_{q^{d_i}}]$
let us denote $\rk_{\F_q}V=\sum_i n_id_i$.

\begin{thm}\label{norm}
One has an isomorphism
$$\FF(F_{V,\chi}(a))\simeq
H\otimes F_{W,\eta}(b)$$ 
where $H$ is a one-dimensional $\ov{\Q}_l$-vector space with
$\Gal(\ov{\F_q}/\F_q)$-action in the following situations:

\noindent
(i) $\rk_{\F_q}V=2$, $W=V$,
$\eta=(\nu\circ\sideset{}{_V}\det)\chi^{-1}$,
where the characters $\nu\in X(\F_q^*)$, $\chi\in X(k^*)$
satisfy
\begin{equation}\label{div1-norm}
D_{1,1}+D_{\nu^{-1},1}=D_{\chi^{-1},V};
\end{equation}
if $d(V)=2$ then we require that $\nu$
is the non-trivial character of order $2$ (so $q$ should be odd);
\begin{equation}\label{ab-norm}
ab=-p(V)^{-1};
\end{equation}
\begin{equation}\label{V1-norm}
H=G(\nu^{-1})\otimes G(\chi)\otimes
(L_{\nu})_{-b}(-m)
\end{equation}
where $m\in\Z$.

\noindent
(ii) $\rk_{\F_q}V=0$, $W=-V$;
$\eta=(\nu\circ\sideset{}{_V}\det)\chi^{-1}$,
where the characters $\nu\in X(\F_q^*)$, $\chi\in X(k^*)$
satisfy
\begin{equation}\label{div2-norm}
D_{1,1}-D_{\nu^{-1},1}=D_{\chi^{-1},V};
\end{equation}
if $d(V)>1$ then we require $\nu=1$;
\begin{equation}\label{a/b-norm}
\frac{a}{b}=p(V);
\end{equation}
\begin{equation}\label{V2-norm}
H=G(\nu)\otimes G(\chi)\otimes
(L_{\nu^{-1}})_{-b}(-m)
\end{equation}
where $m\in\Z$.
\end{thm}

\Pf . For every extension $\F_q\subset \F_{q_1}$ let us
denote $k_{q_1}=k\otimes_{\F_q}\F_{q_1}$. Then
the trace function of the sheaf (\ref{smoothperv}) over $\F_{q_1}$ is given
by
$$f_{V,\chi;q_1}(x)=(-1)^d
\psi(\Tr(a\sideset{}{_{V_{q_1}/\F_{q_1}}}\det(x)))\chi(\Nm_{k_{q_1}/k}(x))
$$
for $x\in k_{q_1}$, where $V_{q_1}=V\otimes_k k_{q_1}$.

An argument similar to that of Lemma \ref{mainlemma} shows that
it suffices to check that for every extension $\F_q\subset\F_{q_1}$ 
and every non-degenerate character $\la$ of $k_{q_1}^*$,  
one has
\begin{equation}
(-q_1)^d(f_{V,\chi;q_1},\la)=\ov{g(\la)} \cdot 
(f_{W,\eta;q_1},\la^{-1})
\end{equation}
and that this number is not zero for at least one 
non-degenerate character $\la$. 
More precisely, we replace the embedding $\G_m^n\ra\A^n$ of Lemma 
\ref{mainlemma} by the embedding $j:\G_mk\ra \A^1k.$
Then we replace the group $\ov{K_0(\G_m^n)}$ from the proof of Lemma
\ref{mainlemma} by the quotient of $K_0(\G_mk)$
by the subgroup $j^*\FF(K_0(Z))$, 
where $Z$ is the complement to the image of $j$.
Then almost the same proof goes through. The only point where one needs
a different argument is in showing that
an element $x\in K_0(\G_mk)$ belongs to
$j^*\FF(K_0(Z))$ if and only if a similar condition holds for all its
trace functions.  For $k$ split over $\F_q$ 
this statement is proven in Lemma \ref{mainlemma}. Thus,
it suffices to check that if
$x\otimes\F_{q_1}\in j^*(\FF(K_0(Z\otimes\F_{q_1})))$ then
$x\in j^*\FF(K_0(Z))$. Equivalently, we have to check that
if an element $y\in K_0(\A^1k)$ satisfies
$$y\otimes\F_{q_1}\in K_0(Z\otimes\F_{q_1})+\FF(K_0(Z\otimes\F_{q_1}))$$
then $y$ itself satisfies the similar condition.
Let us write $y=\sum a_i[F_i]+\sum_j b_j[G_j]$ where
$F_i$ and $G_j$ are simple perverse sheaves, each sheaf $F_i$ satisfies
either $\supp F_i\subset Z$ of $\supp \FF(F_i)\subset Z$, while
each $G_j$ satisfies neither of these conditions. Then
we should have 
$$\sum_j b_j[G_j\otimes\F_{q_1}]\in K_0(Z\otimes\F_{q_1})+
\FF(K_0(Z\otimes\F_{q_1}))$$
which implies $b_j=0$ since the $[G_j\otimes\F_{q_1}]$ 
(and $\FF([G_j\otimes\F_{q_1}])$) 
are linear combinations of simple perverse sheaves
not supported on $Z\otimes\F_{q_1}$.   

The rest of the proof goes similar to that of Theorem \ref{monom}
using Theorem \ref{gaussnorm} and Lemma \ref{intnorm}.
\ed


\begin{thebibliography}{99}
\bibitem{BBD} A.~Beilinson, J.~Bernstein, P.~Deligne, {\it Faisceaux
pervers}, Asterisque 100 (1982), 7--172.
\bibitem{BE1} S.~Bloch, H.~Esnault, {\it Gauss-Manin determinants
for rank one irregular connections on curves}, preprint math.AG/9904088.
\bibitem{BE2} S.~Bloch, H.~Esnault, {\it Gauss-Manin determinant connections
and periods for irregular connections}, preprint math.AG/9912095.
\bibitem{De1} P. Deligne, {\it Applications de la formule des traces
aux sommes trigonom\'etriques}, in {\it Cohomologie Etale (SGA 4 1/2)},
pp. 168--232. Lecture Notes in Math. 569, Springer, 1977.
\bibitem{De2} P. Deligne, {\it Les constantes des \'equations
   fonctionnelles des fonctions $L$}, in {\it Modular functions of one
   variable}, II (Proc. Internat. Summer School, Univ. Antwerp, Antwerp,
   1972), pp. 501--597. Lecture Notes in Math. 349, Springer,
   1973.
\bibitem{DeWeil} P. Deligne, {\it La conjecture de Weil, II},
Publ. Math. IHES 52 (1980), 313--428.
\bibitem{EKP} P. Etingof, D. Kazhdan, A. Polishchuk, {\it When is
the Fourier transform of an elementary function elementary?}, preprint.
\bibitem{GL} O. Gabber, F. Loeser, {\it Faisceaux pervers $l$-adiques sur
un tore}, Duke Math J. 83 (1996), 501--606.
\bibitem{DH} H. Hasse, H. Davenport, {\it Die Nullstelen der
Kongruenzzetafunktionen in gewissen zyklischen Fallen}, J. Reine
Angew. Math. 172 (1934), 151--182.
\bibitem{K1} N. Katz, {\it Gauss sums, Kloosterman sums, and monodromy
groups}, Princeton University Press, 1988.
\bibitem{K2} N. Katz, {\it Exponential sums and differential equations},
Princeton University Press, 1990.
\bibitem{La1} G. Laumon, {\it Majorations de sommes trigonom\'etriques
(d'apr\`es P. Deligne et N. Katz)}, in 
{\it Caracteristique d'Euler-Poincar\'e},  
Ast\'erisque, 83--83 (1981), 221--258.
\bibitem{La2} G. Laumon, {\it Transformation de Fourier, constantes
   d'\'equations fonctionnelles et conjecture de Weil},
   IHES Publ. Math. No. 65 (1987), 131--210.  
\end{thebibliography}
\end{document}